\documentclass[11pt,a4paper]{article}

\usepackage{fancyhdr}
\usepackage{amsmath,mathtools}
\usepackage{physics}
\usepackage{bm}
\usepackage{esint}
\usepackage{amsfonts}
\usepackage{amsthm}
\usepackage{amssymb}
\usepackage{amsbsy}
\usepackage{verbatim}
\usepackage{graphicx}
\usepackage{subfig}
\usepackage{url}
\usepackage{mathrsfs}
\usepackage[hmargin = 1in,vmargin=.75in]{geometry}
\usepackage{color}
\usepackage{amstext}
\usepackage{stmaryrd}
\usepackage{enumerate}
\usepackage{algpseudocode}
\usepackage{algorithm}
\usepackage{pifont}
\usepackage{prettyref}
\usepackage{hyperref}
\hypersetup{
    colorlinks = true,
}
\font\msbm=msbm10

\numberwithin{equation}{section}

\theoremstyle{plain}
\newtheorem{theorem}{Theorem}[section]
\newtheorem{lemma}[theorem]{Lemma}
\newtheorem{corollary}[theorem]{Corollary}

\newtheorem{proposition}[theorem]{Proposition}

\newtheorem{definition}{Definition}[section]

\def\mathbb#1{\hbox{\msbm{#1}}}

\newcommand{\out}[2]{#1#2^{\top}}

\newcommand{\bg}{\boldsymbol{g}}

\newcommand{\bu}{\boldsymbol{u}}

\newcommand{\bx}{\boldsymbol{x}}

\newcommand{\BA}{\boldsymbol{A}}

\newcommand{\BC}{\boldsymbol{C}}

\newcommand{\BE}{\boldsymbol{E}}

\newcommand{\BG}{\boldsymbol{G}}
\newcommand{\BH}{\boldsymbol{H}}
\newcommand{\BI}{\boldsymbol{I}}

\newcommand{\BK}{\boldsymbol{K}}

\newcommand{\BO}{\boldsymbol{O}}

\newcommand{\BQ}{\boldsymbol{Q}}

\newcommand{\BS}{\boldsymbol{S}}

\newcommand{\BV}{\boldsymbol{V}}
\newcommand{\BW}{\boldsymbol{W}}
\newcommand{\BX}{\boldsymbol{X}}
\newcommand{\BY}{\boldsymbol{Y}}
\newcommand{\BZ}{\boldsymbol{Z}}

\newcommand{\BPhi}{\boldsymbol{\Phi}}

\newcommand{\BPi}{\boldsymbol{\Pi}}
\newcommand{\BDelta}{\boldsymbol{\Delta}}

\newcommand{\RR}{\mathbb{R}}

\newcommand{\BLambda}{\boldsymbol{\Lambda}}
\newcommand{\BSigma}{\boldsymbol{\Sigma}}

\newcommand{\PP}{\mathcal{P}}

\newcommand{\I}{\boldsymbol{I}}
\newcommand{\reals}{{\mathbb{R}}}

\newcommand{\lag}{\langle}
\newcommand{\rag}{\rangle}

\newcommand{\lp}{\left(} 
\newcommand{\rp}{\right)} 

\newcommand{\eps}{\epsilon}

\newcommand*\diff{\mathop{}\!\mathrm{d}}

\DeclareMathOperator{\rgrad}{grad}

\DeclareMathOperator{\Proj}{Proj}
\DeclareMathOperator{\Hess}{Hess}

\DeclareMathOperator{\E}{\mathbb{E}}

\DeclareMathOperator{\SNR}{SNR}

\DeclareMathOperator{\St}{St}
\DeclareMathOperator{\sk}{skew}

\DeclareMathOperator{\Od}{O}

\DeclareMathOperator{\eig}{eig}

\DeclareMathOperator*{\argmin}{argmin}

\newcommand{\HBZ}{\widehat{\BZ}}

\newcommand{\BSeig}{\widehat{\BS}^{\rm eig}}
\newcommand{\BSeigloo}{\widehat{\BS}^{\rm eig (i)}}
\newcommand{\BSeigloot}{\widehat{\BS}^{\rm eig (i) \top}}
\newcommand{\BSeigt}{\widehat{\BS}^{{\rm eig}\top}}
\newcommand{\BQeig}{\BQ^{\rm eig}}

\newcommand{\BZeig}{\widehat{\BZ}^{\rm eig}}

\renewcommand{\Pr}{\mathbb{P}}

\begin{document}
\title{\bf Uncertainty Quantification of Spectral Estimator and MLE for Orthogonal Group Synchronization}
\author{Ziliang Samuel Zhong\thanks{Shanghai Frontiers Science Center of Artificial Intelligence and Deep Learning, New York University Shanghai, Shanghai, China. S.L. and Z.S.Z. are (partially) financially supported by the National Key R\&D Program of China, 2021YFA1002800, Shanghai STCSM Rising Star Program 24QA2706200, STCSM General Program 24ZR1455300, National Natural Science Foundation of China 12001372, Shanghai SMEC 0920000112, and NYU Shanghai Boost Fund. Z.S.Z. is also supported by NYU Shanghai Ph.D. fellowship.} \thanks{Center for Data Science, New York University.}, Shuyang Ling$^*$}
\maketitle

\begin{abstract}
Orthogonal group synchronization aims to recover orthogonal group elements from their noisy pairwise measurements. It has found numerous applications including computer vision, imaging science, and community detection. Due to the orthogonal constraints, it is often challenging to find the least squares estimator in presence of noise. 
In the recent years, semidefinite relaxation (SDR) and spectral methods have proven to be powerful tools in recovering the group elements. In particular, under additive Gaussian noise, the SDR exactly produces the maximum likelihood estimator (MLE), and both MLE and spectral methods are able to achieve near-optimal statistical error. 
In this work, we take one step further to quantify the uncertainty of the MLE and spectral estimators by considering their distributions. 
By leveraging the orthogonality constraints in the likelihood function, 
we obtain a second-order expansion of the MLE and spectral estimator with the leading terms as an anti-symmetric Gaussian random matrix that is on the tangent space of the orthogonal matrix. This also implies state-of-the-art min-max risk bounds and a confidence region of each group element as a by-product. Our works provide a general theoretical framework that is potentially useful to find an approximate distribution of the estimators arising from many statistical inference problems with manifold constraints. The numerical experiments confirm our theoretical contribution.

\end{abstract}





\section{Introduction}

Group synchronization aims to recover the unknown group elements $\bg_i\in{\cal G}$ from their pairwise measurements:
\begin{equation}
\label{eq:groupsync}
	\bg_{ij} = \bg_{i}\bg_{j}^{-1}+ \text{noise}.
\end{equation}
Depending on the specific group ${\cal G}$, the group synchronization problem has found numerous applications in practice, e.g. $\mathbb{Z}_2$-synchronization for community detection~\cite{A17,B18}, angular synchronization for clock synchronization~\cite{BBS17,S11,ZB18}, statistical ranking~\cite{C16,DCT21}, point cloud registration~\cite{CKS15,L23a}, and cryo-EM microscopy~\cite{S18}.

In this work, we will focus on the orthogonal group synchronization, i.e., ${\cal G} = \Od(d)$ where
\[
\Od(d) := \{\BO \in \reals^{d \times d}: \BO\BO^{\top} = \BO^{\top}\BO  = \BI_{d} \}
\]
denotes the set of $d\times d$ orthogonal matrices. More precisely, for every $(i, j) \in [n]^2$, we observe the noisy pairwise measurements between two group elements: 
$\BC_{ij} = \BZ_{i}\BZ_{j}^{\top} + \BDelta_{ij}$. The compact form is given by $\BC = \BZ\BZ^{\top} + \BDelta$
where $\BZ^{\top} = [\BZ_1^{\top},\cdots,\BZ_n^{\top}] \in \Od(d)^{\otimes n}$ is the stacking of $\BZ_i\in\Od(d)$, $i \in [n]$.
 It is particularly interesting as it includes $\mathbb{Z}_2$-synchronization, angular synchronization, and point cloud registration as spectral examples. 
 
\noindent{\bf Least squares estimator}~~ Given this inverse problem, a common approach is to find the least squares (LS) estimator via minimizing the following optimization problem:
\begin{equation}
\label{eq:ls}
	\min_{\BS_{1}, \ldots, \BS_{n} \in \Od(d)} \sum_{i<j} \norm{\BS_{i}\BS_{j}^{\top} - \BC_{ij}}^{2}_{F} \iff \min_{\BS_{1}, \ldots, \BS_{n} \in \Od(d)} -\lag \BC, \BS\BS^{\top} \rag.
	\tag{LS}
\end{equation}
where $\BS^{\top} = [\BS_{1}^{\top}, \ldots, \BS_{n}^{\top}] \in \Od(d)^{\otimes n}$. Finding the exact solution to~\eqref{eq:ls} can be hard since even when $d=1$,~\eqref{eq:ls} is equivalent to minimizing a quadratic form over hypercube $\{\pm 1\}^n$. Therefore, in the worst case, it can be potentially NP-hard to find the global minimizer~\cite{GW95}. 

\noindent{\bf Semidefinite relaxation (SDR)~~}
Semidefinite relaxation is a powerful tool that often provides an accurate approximation of many nonconvex problems. In particular,~\eqref{eq:ls} has a natural SDR inspired by~\cite{GW95}. 
Let $\BX = \out{\BS}{\BS} \in \reals^{nd \times nd}$ with $\BX_{ij} = \out{\BS_{i}}{\BS_{j}}$. Note that $\BX_{ii} = \BI_{d}$ and $\BX \succeq 0$ hold for any $\BS_{i} \in \Od(d)$, and then it gives an SDR of~\eqref{eq:ls} 
\begin{equation}\label{eq:sdp}
	\min_{\BX \in \reals^{nd \times nd}} - \lag \BC,  \BX\rag \qquad \text{s.t} \qquad \BX_{ii} = \BI_{d}, \quad \BX \succeq 0.
\tag{SDR}
\end{equation}

\noindent{\bf Spectral relaxation}~~ Another useful approximation to~\eqref{eq:ls} is the spectral relaxation. Note that $\BS \in \Od(d)^{\otimes n}$ satisfies $\BS^{\top}\BS = n\BI_{d}$. The spectral method simply replaces the constraint $\BS \in \Od(d)^{\otimes n}$ in~\eqref{eq:ls} with $\BS^{\top}\BS = n\BI_{d}$, which reads
\begin{equation} \label{eq:spectral}\tag{SPR}
	\begin{aligned}
		\min &\quad -\lag \BC, \BS\BS^{\top}\rag\qquad {\rm s.t.} \qquad \BS^{\top}\BS = n\BI_{d}.
	\end{aligned}
\end{equation}
We denote $\widehat{\BS}^{\rm eig} = [\BS_{1}^{{\rm eig} \top}, \ldots, \BS_{n}^{{\rm eig} \top}]^{\top}$ as its global minimizer, i.e., the top $d$ eigenvectors of $\BC$. Since each $\widehat{\BS}_i^{\eig}$ is not necessarily orthogonal, a rounding procedure $\BZeig := \PP_{n}(\widehat{\BZ}^{\rm eig}) \in \Od(d)^{\otimes n}$ is needed to give the final spectral estimator to $\BZ$ where this rounding operator ${\cal P}_n(\cdot)$ will be formally introduced in~\eqref{def:Pn}.

With the above two relaxations, one natural question is how~\eqref{eq:sdp} and~\eqref{eq:spectral} approximates the least squares estimator~\eqref{eq:ls}.
As a benchmark model, orthogonal group synchronization under the additive Gaussian noise is widely considered:
\begin{equation}\label{eq:model}
\BC_{ij} = \BZ_{i} \BZ_{j}^{\top} +\sigma\BW_{ij}~~\Longleftrightarrow~~\BDelta = \sigma\BW
\end{equation}
where  $\BW \in \reals^{nd \times nd}$ is a symmetric random matrix whose diagonal blocks are zero, i.e., $\BW_{ii} = 0$ and  and $\BW_{ij} \in \reals^{d\times d}$ is the $(i,j)$-th block of $\BW$ with independent $\mathcal{N}(0,1)$ entries satisfying $\BW_{ij}^{\top} = \BW_{ji}$. 
It is important to note that under Gaussian noise,~\eqref{eq:ls} is equivalent to the maximum likelihood estimator (MLE).


\noindent{\bf Tightness of SDR}~~
For $d\geq 2$, the SDR does not yield the ground truth, unlike the $d=1$ case~\cite{A17,B18}.
It is surprising that the SDR is tight~\cite{BBS17,LYS17,ZB18,L22} if the noise is small compared with the signal strength, i.e., it recovers the MLE with high probability  when
\begin{align*} 
	\sigma \leq c_0\sqrt{\frac{n}{d}} \cdot\frac{1}{\sqrt{d} + \sqrt{\log n}}
\end{align*}
for some absolute positive constant $c_0 > 0.$  In other words, the tightness of the SDR means the estimator $\widehat{\BX}$ is exactly rank-$d$ and can be factorized into $\widehat{\BX} = \widehat{\BZ}\widehat{\BZ}^{\top}$ for some $\widehat{\BZ}\in\Od(d)^{\otimes n}$, and $\widehat{\BZ}$ is the MLE. This is information-theoretically near-optimal modulo an extra factor $\sqrt{d} + \sqrt{\log n}$. The proof follows from the convergence analysis of the generalized power method~\cite{LYS17,LYS23,L22,ZB18} whose justification is combined with the leave-one-out technique~\cite{AFWZ20}.

\noindent{\bf Near-optimal spectral estimator and SDP estimator}~~
For both the spectral estimator $\widehat{\BZ}^{\eig}$ and SDP estimator $\widehat{\BZ}$, their deviations from the ground truth $\BZ$ are controlled by
\[
\min_{\BO\in\Od(d)} \|\widehat{\BZ} - \BZ\BO\|_F\lesssim \sigma\sqrt{d} (\sqrt{d} + \sqrt{\log n}),~~~\min_{\BO\in\Od(d)} \|\widehat{\BZ}^{\eig} - \BZ\BO\|_F\lesssim \sigma\sqrt{d} (\sqrt{d} + \sqrt{\log n}).
\]
The error bound between the MLE $\widehat{\BZ}$ and the ground truth $\BZ$ is obtained in~\cite{L22} and the counterpart for angular synchronization is given in~\cite{ZB18}. 
Moreover, it has been shown in~\cite{GZ22,GZ23,Z22} that under min-max risk, the error bound above is near-optimal.

Despite those aforementioned progresses, a more precise statistical characterization of the MLE and spectral estimator is lacking. From a statistical perspective, obtaining the approximate distribution of the error is crucial for constructing a confidence interval of each group element $\BZ_i$. 
This observation motivates this work which leads to the following question:
\begin{equation}
\text{\em What is the approximate distribution of MLE and spectral estimator?}
\end{equation}
In other words, we are interested in obtaining the precise statistical property of the minimizer of the MLE and spectral estimator. The challenges in quantifying the uncertainty of $\widehat{\BZ}$ and $\widehat{\BZ}^{\eig}$ lie on the fact that the likelihood function is nonconvex due to the orthogonality constraints, unlike most of the existing literatures regarding the uncertainty quantification in many other high-dimensional inference problems. The inference under orthogonal constraints is mainly known for eigenvector and eigenvalue inference~\cite{BBP05,BDWW22}. In the past few years, there are also many works focusing on the entrywise analysis of the eigenvectors in community detection and synchronization~\cite{A17,AFWZ20,L22b} that provide the first-order error expansion while quantifying the uncertainty of the estimators is essentially a second-order error expansion.

The precise characterization of the statistical estimator is also known as the uncertainty quantification. The uncertainty quantification has recently gained much attention in a series of applications such as statistical ranking~\cite{GSZ23} and matrix completion~\cite{CFM+19}. 
In high-dimensional linear regression, there have quite a few important works in constructing the confidence intervals of the estimators, and derived the asymptotic distribution for the Lasso estimator under high-dimensional model, e.g.~\cite{ZZ14,CG17,JM14,VSV17}. For the generalized linear model,~\cite{SC19} characterizes the distribution of likelihood ratio statistics and the asymptotic bias and variance of the MLE in the high-dimensional logistic regression.
In graphical model,~\cite{JV15,RSZ15} proposed methodology to estimate the sparse inverse covariance estimation in the high-dimensional settings, and obtained the limiting distribution of the estimators.

Therefore, the settings proposed in our work have little overlaps with the existing literature. Our work features the following highlights: we obtain a precise asymptotic distribution of both estimators, which extends state-of-the-art first-order entrywise analysis of the MLE and spectral estimators. The results imply that both are near-optimal under the min-max risk as a by-product. Moreover, we show that the MLE achieves a slightly better second-order error than the spectral estimator does, although their first-order terms are the same.
While the current work only considers the additive Gaussian noise, it has great potentials to extend to more general scenarios such as sub-gaussian noise.
Our works provide a general theoretical framework that is potentially useful to find the approximate distribution of the estimators arising from many statistical inference problems with manifold constraints.

\subsection{Notation}
 We denote boldface $\bx$ and $\BX$ as a vector and matrix respectively, and $\bx^{\top}$ and $\BX^{\top}$ are their corresponding transpose. The $\I_n$ stands for the $n\times n$ identity matrix. 
For two matrices $\BX$ and $\BY$ of the same size, $\lag \BX,\BY\rag = \sum_{i,j}X_{ij}Y_{ij}$ is their inner product. For any matrix $\BX$, the Frobenius, operator, and nuclear norms are denoted by $\|\BX\|_F$, $\|\BX\|$, and $\|\BX\|_*$ respectively. For two nonnegative functions $f(x)$ and $g(x)$, we say $f(x)\lesssim g(x)$ (and $f(x) \gtrsim g(x)$) if there exists a universal constant $C> 0$ such that $f(x) \leq Cg(x)$ (and $f(x)\geq Cg(x)$). We say a matrix $\BX = O(T)$ if $\norm{\BX}_{F} \le c \cdot T$ for some constant $c>0$.

\subsection{Organization}
Section~\ref{s:main} gives the main results of our work, and Section~\ref{s:numerics} provides numerical experiments that complement our theory. The proof is deferred to the Appendix.

\section{Main theorems}\label{s:main}
Before we proceed to the main results, we introduce a few important notations.
Throughout our presentation, we will often use the quantity $\SNR$, i.e., the signal-to-noise ratio, which is defined as the ratio between the operator norms of the signal and the noise part in~\eqref{eq:model}:
\begin{equation}
	\SNR := \sigma^{-1}\sqrt{\frac{n}{d}} \asymp \frac{\norm{\BZ \BZ^{\top}}}{\sigma\norm{\BW}}
\end{equation}
where $\|\BW\| = (2+o(1))\sqrt{nd}$ and $\|\BZ\BZ^{\top}\| = n.$

Note that for any estimator $\BS$ of $\BZ\in\Od(d)^{\otimes n}$, we note that $\BS\BO$ satisfies $\BS\BO\BO^{\top}\BS^{\top} = \BS\BS^{\top}$ for any orthogonal matrix $\BO$. Therefore, to resolve this ambiguity, we define the following distance function for any $\BS$ and $\BZ \in \reals^{nd \times d}$:
\[
d_F(\BS,\BZ) : = \min_{\BO\in\Od(d)}\|\BS - \BZ\BO\|_F.
\] The minimum is obtained via the polar retraction, which is defined as $\mathcal{P}: \reals^{d \times d} \mapsto \Od(d)$,
\begin{equation}\label{def:polar0}
\mathcal{P}(\BY) \in \argmin_{\BO\in\Od(d)}\|\BO - \BY\|_F^2, 
\end{equation}
where $\BY \in \reals^{d\times d}$. 
In fact, $\mathcal{P}(\BY)$ equals the product of the left and right singular vectors  of $\BY$ and ${\cal P}(\cdot)$ is not unique if $\BY$ is rank-deficient. If $\BY$ is full-rank, then
\begin{equation}\label{def:polar}
{\cal P}(\BY) := (\BY\BY^{\top})^{-1/2}\BY = \BY (\BY^{\top}\BY)^{-1/2}.
\end{equation}
We also define $\mathcal{P}_{n}: \reals^{nd \times d} \mapsto \Od(d)^{\otimes n}$ which applies ${\cal P}$ to each $d\times d$ block:
\begin{equation}\label{def:Pn}
\mathcal{P}_{n}([\BY^{\top}_{1},\ldots,\BY^{\top}_{n}]^{\top}) := [\mathcal{P}(\BY_{1})^{\top},\ldots,\mathcal{P}(\BY_{n})^{\top}]^{\top}.
\end{equation}

\subsection{Distributional guarantees}\label{ss: distributional guarantees}
\subsubsection{Distribution of the MLE}

Our first theorem gives a block-wise approximate distribution of each $\widehat{\BZ}_i$, $1\leq i\leq n$, where $\widehat{\BZ}_i$ is the MLE of $\BZ_i.$ 
\begin{theorem}[\bf Distribution of the MLE] \label{thm:main}
	 Let $\widehat{\BZ} = [\widehat{\BZ}^{\top}_{1}, \ldots, \widehat{\BZ}^{\top}_{n}]^{\top}$ be the global minimizer~\eqref{eq:ls} and $d\geq 2$. Suppose 
	 \begin{align} \label{eqn: snr}
	 	\SNR \gtrsim \sqrt{d} + \sqrt{\log n}, 
	 \end{align}
	 then for $1\le i \le n$, each $\widehat{\BZ}_{i}$ has the following decomposition:
	\begin{equation}\label{eq:main}
		\widehat{\BZ}_{i}  = \left(\I_d + \frac{\sigma\sqrt{n-1}}{2n} (\BG_i - \BG_i^{\top}) \right)\BZ_i\BQ + \BE_i			\end{equation}
where $\BG_{i}$ is a $d\times d$ Gaussian random matrix with i.i.d. $\mathcal{N}(0,1)$ entries and $\BQ = \mathcal{P}(\BZ^{\top}\widehat{\BZ})$, and $\BE_i$ satisfies with probability at least $1-n^{-2}$ that
\[
\|\BE_i\|_F \lesssim \SNR^{-2}(\sqrt{d} + \sqrt{\log n})^2.
\]

\end{theorem}

Here are a few remarks regarding this main theorem. Note that the right hand side of~\eqref{eq:main} involves an unknown orthogonal matrix $\BQ$ which depends on $\BZ$. This is the inherent ambiguity in all group synchronization problems because the recovery of the group element is only unique up to a global group action due to the pairwise observations. One immediate application of this result is to derive the distribution of $\widehat{\BZ}_i\widehat{\BZ}_j^{\top}$ and then obtain a confidence region of $\BZ_i\BZ_j^{\top}$. The approximate distribution of $\widehat{\BZ}_i\widehat{\BZ}_j^{\top}$  will not be affected by the unknown matrix $\BQ$, as will be discussed in details later.

In addition, the reason why we only consider $d\geq 2$ is because for $d=1$, the SDR will exactly recover the ground truth. Unlike $d=1$, the global minimizer given by the SDR for $d\geq 2$ is typically not equal to the ground truth. Therefore, it makes more sense to consider $d\geq 2.$ 
Our result proves the asymptotic normality of the MLE under the information-theoretically near-optimal $\SNR$. 
The work~\cite{L22} provided an error bound between the MLE $\widehat{\BZ}$ and the ground truth $\BZ$ by
\[
\min_{\BO\in\Od(d)}\max_{1\leq i\leq n}\| \widehat{\BZ}_i - \BZ_i\BO \|_F \lesssim \SNR^{-1} (\sqrt{d} + \sqrt{\log n}).
\]
Moreover, it was shown in~\cite{GZ23} by Gao and Zhang that the MLE is optimal in terms of a min-max risk bound under $d_F(\cdot,\cdot)$:
\[
\inf_{\BS\in\Od(d)^{\otimes n}} \sup_{\BZ \in  \Od(d)^{\otimes n}}\E_{\BZ}  d_F^2(\BS, \BZ) \geq  \frac{(1-\delta)\sigma^2 d(d-1)}{2}
\]
where $\BS$ is an estimator of $\BZ$. A similar result also holds for spectral estimator~\cite{Z22} under the same setting.
If we look into~\eqref{eq:main}, it is straightforward to obtain
\begin{align*}
\E\|\widehat{\BZ}_i - \BZ_i\BQ\|_F^2 & = \frac{\sigma^2(n-1)\cdot 2d(d-1)}{4n^2} + \frac{\sigma(\sqrt{d} + \sqrt{\log n})}{\sqrt{n}}\cdot  O \lp \SNR^{-2}(\sqrt{d} + \sqrt{\log n})^2 \rp  \\
& = \frac{\sigma^2d(d-1)}{2n} + O \lp \frac{(\SNR^{-1}(\sqrt{d} + \sqrt{\log n}))^3}{\sqrt{d}} \rp 
\end{align*}
where $\E \|\BG_i-\BG_i^{\top}\|_F^2= 2d(d-1)$, $\|\BG_i - \BG_i^{\top}\| \lesssim \sqrt{d} + \sqrt{\log n}$, and $\sigma^2 d = \SNR^{-2}n$.
This implies
\begin{align*}
\E\|\widehat{\BZ} - \BZ\BQ\|_F^2 
& = \frac{\sigma^2 d(d-1)}{2} \left(  1 +  O\left(\SNR^{-1}\left(1 + \sqrt{\frac{\log n}{d}}\right)^3\right)\right).
\end{align*}
As a result, our theorem provides a near-optimal characterization of the distribution of each $\widehat{\BZ}_i$ and~\eqref{eq:main} implies the near-optimal min-max risk bound under $d_F(\cdot,\cdot).$

Finally, we give an intuitive interpretation that naturally explains the second-order expansion in~\eqref{eq:main}. Note that $\BZ_{i}\BQ$ is on the manifold $\Od(d)$ and the tangent space at $\BZ_i\BQ$ for any $\BQ\in\Od(d)$ with $d\geq 2$ is given by
\[
\{\BV_{i}\in\RR^{d\times d}: \BV_{i}\BQ^{\top} \BZ^{\top}_{i} + \BZ_{i}\BQ \BV^{\top}_{i} =0 \}.
\]
The local expansion of $\widehat{\BZ}_i$ in Theorem~\ref{thm:main} can be expressed into
\[
\widehat{\BZ}_{i} = \BZ_{i}\BQ + \underbrace{\frac{\sigma\sqrt{n-1}}{2n} (\BG_i - \BG_i^{\top})\BZ_i\BQ}_{\text{tangent vector}} +\text{second order error},
\]
i.e., the leading error term is on the tangent space of $\Od(d)$ at $\BZ_i\BQ.$
Therefore, our result provides a characterization of the MLE that is subject to the orthogonality constraints via tangent space at the ground truth.

\subsubsection{Distribution of the spectral estimator}
We also derive the approximate distribution for the spectral estimator. Let $\widehat{\BS}^{\eig}$ be the global minimizer to~\eqref{eq:spectral}, and then the spectral estimator is given by $\widehat{\BZ}_i^{\eig} = {\cal P}(\BS_i^{\eig})$
where ${\cal P}(\cdot)$ is the polar retraction~\eqref{def:polar0}.
Similar to the MLE, a block-wise error bound has been obtained in~\cite{L22b}:
\[
\min_{\BO\in\Od(d)} \max_{1\leq i\leq n} \|\widehat{\BZ}_i^{\eig} - \BZ_i\BO\|_F\lesssim \SNR^{-1}(\sqrt{d} +\sqrt{\log n})
\]
provided that $\SNR\gtrsim \sqrt{d} +\sqrt{\log n}.$ The following theorem gives a precise distribution of $\widehat{\BZ}_i^{\eig}$, i.e., a blockwise second-order expansion of $\widehat{\BZ}_i^{\eig}$.

\begin{theorem}[\bf Distribution of the spectral estimator] \label{thm:spectral_main}
	 Let $\widehat{\BS}^{\rm eig} = [\BS_{1}^{{\rm eig} \top}, \ldots, \BS_{n}^{{\rm eig} \top}]^{\top}$ be the solution of~\eqref{eq:spectral} and $d\geq 2$. Suppose 
\[
\SNR \gtrsim \sqrt{\kappa} (\sqrt{d} + \sqrt{\log n}), 
\]
	 then for $1\le i \le n$, each $\widehat{\BZ}^{\rm eig}_{i} := \PP(\BSeig_{i}) $ satisfies the following decomposition
	\begin{equation}\label{eq:main_spectral}
		\widehat{\BZ}^{\rm eig}_{i} =\left(\I_d +  \frac{\sigma\sqrt{n-1}}{2n}(\BG_i - \BG_i^{\top})\right)\BZ_i\BQ^{\rm eig} + \BE_i^{\eig}
	\end{equation}
	where 
$\BG_{i}$ is a $d\times d$ Gaussian random matrix with i.i.d. $\mathcal{N}(0,1)$ entries, $\BQ^{\rm eig} = \mathcal{P}(\BZ^{\top}\widehat{\BZ}^{\rm eig}),$
and $\BE_i$ satisfies with high probability that
\[
\|\BE_i^{\eig}\|_F \lesssim  \kappa^{1/2} \SNR^{-2}(\sqrt{d}+\sqrt{\log n})^2,~~1\leq i\leq n,
\]
where
$\kappa$ is the condition number of $\BS_i^{\eig\top}\BS_i^{\eig}$:
\[
\kappa = \max_{1\leq i\leq n} \frac{\lambda_{\max}(\BS_i^{\eig\top}\BS_i^{\eig})}{\lambda_{\min}(\BS_i^{\eig\top}\BS_i^{\eig})} \leq \frac{1 + c'\SNR^{-1}}{1 - c'\SNR^{-1}}
\]
for some constant $c'>0.$
\end{theorem}
The remarks to Theorem~\ref{thm:main} also apply to the spectral case. 
Compared to~\eqref{eq:main}, the error term of~\eqref{eq:main_spectral} contains an extra scalar $\kappa \geq 1$, the condition number of $\BSeigt_{i}\BSeig_{i}$. In other words, $\SNR$ is sufficiently large, the second-order expansion of MLE and spectral estimators are nearly the same. However, for $\SNR$ moderately large, then the second-order error of spectral estimator is strictly larger than that of the MLE. Similar to the MLE, the spectral estimator also achieves the min-max optimal risk.

We conclude this section by pointing out a few future directions. Our work only focuses on the orthogonal group synchronization under Gaussian noise. There are a few future directions that are worth exploring. 
\begin{itemize}
\item Extension to other sub-gaussian noise for least squares estimator: To extend the result for other type of noise, we need to first show that the SDR is tight and recovers the LS estimator if the SNR is large; once we have that, we can approximate the least squares estimator via one-step Newton iteration as employed in this work. Then it suffices to estimate the distribution of one-step Newton iteration via the central limit theorem.

\item Extension to other sub-gaussian noise for spectral estimator: for spectral estimator, a similar argument will also work for $\Od(d)$-synchronization with uniform corruption model. The main difference is about establishing the central limit theorem for a sum of random matrices. 

\item The current framework is very versatile and could provide a simple way to derive the approximate distribution of the estimators in other statistical inference problems with manifold constraints, including the inferences about eigenvector/eigenvalue that arise from high-dimensional statistics.
\end{itemize}

\subsection{Confidence region}
{

In this section, we illustrate how to construct confidence regions using the second-order expansions in Theorem~\ref{thm:main} and~\ref{thm:spectral_main}. To concisely present a unified approach to establish confidence regions for both MLE and spectral estimator, we use the notation $\widehat{\BS} \in \{\widehat{\BZ},\BZeig\}$ to denote either the MLE or spectral estimator. Recall that  $\widehat{\BS}$ satisfies
\begin{equation}\label{eq: confidence region mle}
    \BZ_{i}\BQ \widehat{\BS}_{i}^{\top} = \BI_{d} + \frac{\sigma \sqrt{n-1}}{2n} \lp \BG_{i}^{\top} - \BG_{i} \rp + \BZ_{i}\BQ \BE_{i}^{\top},~~\BQ = \PP(\BZ^{\top}\widehat{\BS}),~~1\leq i\leq n,
\end{equation}
with probability at least $1 -O( n^{-2})$ where $\norm{\BE_{i}}_{F} \lesssim \SNR^{-2} \lp \sqrt{d} + \sqrt{\log n} \rp^2$. As discussed in Section~\ref{ss: distributional guarantees}, the anti-symmetric Gaussian term dominates~\eqref{eq: confidence region mle}, i.e., $  \SNR^{-1}(\sqrt{d}+\sqrt{\log n}) \asymp \frac{\sigma \sqrt{n-1}}{2n} \norm{ \BG_{i}^{\top} - \BG_{i}}_{F} \geq \norm{\BE_{i}}_{F} \asymp  \SNR^{-2}(\sqrt{d}+\sqrt{\log n})^2$. 
Now we proceed to present the construction of the confidence region of $\BZ_i.$

\noindent\textbf{Estimation of $\sigma^2$.} We estimate $\sigma^2$ by using the residue:
\begin{equation}\label{eq: variance estimation}
    \widehat{\sigma}^2 = \frac{1}{n(n-1)d^2} \|\widehat{\BS}\widehat{\BS}^{\top} - \BC\|^2_F.
\end{equation}

\noindent\textbf{Confidence region.} 
Motivated by the approximate normality of~\eqref{eq: confidence region mle}, we propose the following approximate $(1-\alpha)$-confidence regions of $\BZ_i\BQ$ that is given by ${\rm CR}_{1-\alpha}\widehat{\BS}_i$, i.e., 
\[
{\rm CR}_{1-\alpha}\widehat{\BS}_i : = \left\{ \PP(\I_d + \BA)\widehat{\BS}_i: \|\BA\|_F \leq  n^{-1/2}\widehat{\sigma} \cdot c_{1-\alpha}  \right\}
\]
where
\begin{equation}\label{eq: CR mle}
\begin{aligned}
{\rm CR}_{1-\alpha} & : = \left\{ \PP(\I_d + \BA): \|\BA\|_F \leq  n^{-1/2}\widehat{\sigma} \cdot c_{1-\alpha}  \right\},\\
c_{1-\alpha} & := (\chi^2_{d(d-1)/2, 1-\alpha})^{1/2}  + C_0 \widehat{\sigma} n^{-1/2}d (\sqrt{d} + \sqrt{\log n})^2,
\end{aligned}
\end{equation}
for some constant $C_0$ and $\chi^2_{d(d-1)/2, 1-\alpha}$ is the $(1-\alpha)$-quantile of the chi-squared distribution with $d(d-1)/2$ degrees of freedom. In particular, if the noise is small, i.e., $\sigma\ll \sqrt{n}/(\sqrt{d} + \sqrt{\log n})^2$, then $c_{1-\alpha} \approx (\chi^2_{d(d-1)/2, 1-\alpha})^{1/2}$ provides an accurate approximation. 
The following corollary establishes the validity of the confidence region in~\eqref{eq: CR mle}.
\begin{corollary}[Coverage probability of the confidence region] \label{cor: CR mle} 
For $\widehat{\BS} \in \{\widehat{\BZ},\BZeig\}$, under the conditions of Theorem~\ref{thm:main} (or~\ref{thm:spectral_main}), we have:
\[
\mathbb{P} \lp \BZ_{i}\BQ \in  {\rm CR}_{1-\alpha} \widehat{\BS}_i \rp \ge  1-\alpha + o(n^{-1})
\]
for each fixed $1\leq i\leq n$. For any $i\neq j$, a $(1-\alpha)$-confidence region of $\BZ_i\BZ_j^{\top}$ is given by
\[
{\rm CR}_{1-\alpha/2} \widehat{\BS}_i\widehat{\BS}_j^{\top}{\rm CR}_{1-\alpha/2} := \{ \PP(\I_d + \BA_i)\widehat{\BS}_i\widehat{\BS}_j^{\top}\PP(\I_d + \BA_j): \|\BA_i\|_F,\|\BA_j\|_F\leq n^{-1/2}\widehat{\sigma}\cdot c_{1-\alpha/2} \}.
\]
\end{corollary}
Note that it is impossible to identify $\BQ$ because of the unknown $\BZ$. This is also due to the global ambiguity in the group synchronization problem: the recovery is unique up to a global group action. Therefore, the construction of the confidence region for each $\BZ_i$ will involve a common orthogonal matrix $\BQ.$  A confidence region of $\BZ_i\BZ_j^{\top}$ follows directly from that of $\BZ_i\BQ.$
The proof of Corollary~\ref{cor: CR mle} is straightforward.


\begin{proof}[\bf Proof of Corollary~\ref{cor: CR mle}]
Consider the main results in Theorem~\ref{thm:main} and~\ref{thm:spectral_main}: it holds with probability at least $1-O(n^{-2})$ that
\begin{align*}
\BZ_{i}\BQ \widehat{\BS}_{i}^{\top} & = \BI_{d} + \frac{\sigma \sqrt{n-1}}{2n} \lp \BG_{i}^{\top} - \BG_{i} \rp + \BZ_{i}\BQ \BE_{i}^{\top} 
\end{align*}
where $\SNR^{-1} = \sigma \sqrt{d/n}$ and 
\[
\|\BZ_i\BQ\BE_i^{\top}\|_F \leq C_0 \SNR^{-2}(\sqrt{d} + \sqrt{\log n})^2 = C_0 \sigma^2 n^{-1} d (\sqrt{d} + \sqrt{\log n})^2
\] 
for some $C_0>0.$

Consider
\begin{align*}
& \Pr\left( \| \BZ_i \BQ\widehat{\BS}_i^{\top} - \I_d\|_F \geq \widehat{\sigma}n^{-1/2} c\right) \\
& \leq \Pr\left( \frac{\sigma \sqrt{n-1}}{2n}  \left\| \BG_{i}^{\top} - \BG_{i} \right\| + \left\| \BZ_{i}\BQ \BE_{i}^{\top} \right\|_F \geq \widehat{\sigma}n^{-1/2} c  \right) \\
& \leq \Pr\left( \frac{1}{2}\|\BG_i - \BG_i^{\top}\|_F \geq  \frac{\sigma}{\widehat{\sigma}}\cdot c  - \frac{n}{\sigma\sqrt{n-1}}  \left\| \BE_{i} \right\|_F    \right) \\
& \leq \Pr\left( \frac{1}{2}\|\BG_i - \BG_i^{\top}\|_F \geq (1-n^{-1}\lambda )\left(c  -  C_0 \widehat{\sigma} n^{-1/2}d(\sqrt{d} + \sqrt{\log n})^2\right)  \right) 
\end{align*}
where $\BG_i - \BG_i^{\top}$ is anti-symmetric and each entry is $\mathcal{N}(0,2)$,  and 
\[
\frac{1}{4}\|\BG_i - \BG_i^{\top}\|_F^2 = \chi^2_{d(d-1)/2},~~~~\left|\frac{\sigma}{\widehat{\sigma}} - 1\right| \lesssim n^{-1}(1 + d^{-1}\sqrt{\log n}) =: n^{-1}\lambda
\]
which follows from Proposition~\ref{prop: estimation error of sigma2}.
To ensure a coverage probability of $1-\alpha$, we set $c:=c_{1-\alpha}$ such that
\begin{align*}
c_{1-\alpha} 
& = (\chi^2_{d(d-1)/2, 1-\alpha})^{1/2} + C_0 \widehat{\sigma} n^{-1/2}d(\sqrt{d} + \sqrt{\log n})^2 + O(n^{-1}).
\end{align*}

Therefore, with probability at least $1-\alpha$, the set
$\{\I_d+\BA: \|\BA\|_F\leq \widehat{\sigma}n^{-1/2}c_{1-\alpha}\}$
covers $\BZ_i\BQ\widehat{\BS}_i^{\top}$. 
Note that for any $t>0$, it holds that
\[
\left\{ \PP(\I_d + \BA): \|\BA\|_F\leq t\right\} \subseteq  \left\{ \I_d + \BA: \|\BA\|_F\leq t  \right\},
\]
which follows from $\|\PP(\I_d + \BA) - \I_d\|_F \leq \|\BA\|_F.$
In addition, since $\BZ_i\BQ\widehat{\BS}_i^{\top}$ is orthogonal, we have 
\[
\{ \PP(\I_d+\BA): \|\BA\|_F\leq \widehat{\sigma}n^{-1/2}c_{1-\alpha}\}
\]
includes $\BZ_i\BQ\widehat{\BS}_i^{\top}$ with probability $1-\alpha.$
\end{proof}
}

\subsection{Proof sketch of Theorem~\ref{thm:main}: distribution of the MLE} \label{ss:proof sketch mle}

In this section, we provide an overview of the key proof steps. The detailed justifications are provided in Section~\ref{s:riemann} and~\ref{ss:propkey}. Since the MLE does not have a closed-form solution, directly analyzing its  distribution is challenging. To mitigate this issue, we approximate the MLE by the one-step Newton iteration at $\BZ_{i}$ applied to the local likelihood function, defined as:
\begin{equation}
\label{eqn: local MLE}
	  f(\BPhi) := - \sum_{j=1}^n \left \lag \BC_{ij}\widehat{\BZ}_{j}, \BPhi  \right \rag,~~~~\BPhi\in\Od(d)
\end{equation}
where $\sum_{j=1}^n \BC_{ij}\widehat{\BZ}_j = \BZ_i\BZ^{\top}\widehat{\BZ} + \sigma\BW_i^{\top}\widehat{\BZ}.$
Let $\widehat{\BZ} = [\widehat{\BZ}^{\top}_{1}, \ldots, \widehat{\BZ}^{\top}_{n}]^{\top}$ be the global minimizer of~\eqref{eq:ls} with $\widehat{\BZ}_i$ being the $i$-th block of $\widehat{\BZ}$ and it is straightforward to see that
\[
\widehat{\BZ}_i : =\underset{\BPhi\in \Od(d)}{\argmin}~f(\BPhi).
\]

We denote $\widetilde{\BZ}_{i}$ as the Newton approximation to $\widehat{\BZ}_{i}$ at $\BZ_{i}$. Since the objective function $f(\BPhi)$ is defined on the orthogonal matrices, we essentially need some basic tools from Riemannian geometry.
The derivation of the Newton step $\widetilde{\BV}_{i}$ relies on computing the Riemannian gradient and Hessian of $f(\BPhi)$. 
For $f(\BPhi)$ defined above, its Riemannian gradient and Hessian are given by:
	\begin{align*}
	\rgrad f(\BPhi)
	&= -\sum_{j=1}^n \BC_{ij}\widehat{\BZ}_{j} + \frac{1}{2}\sum_{j=1}^n \lp \BC_{ij}\widehat{\BZ}_{j}\BPhi^{\top} +\BPhi\widehat{\BZ}^{\top}_{j}\BC_{ji}\rp \BPhi,\\
	\Hess f(\BPhi)[\dot{\BPhi}, \dot{\BPhi}] &= \left \lag \frac{1}{2}\sum_{j=1}^n \lp \BC_{ij}\widehat{\BZ}_{j}\BPhi^{\top} + \BPhi\widehat{\BZ}^{\top}_{j}\BC_{ji}\rp , \dot{\BPhi}\dot{\BPhi}^{\top} \right \rag,
\end{align*}
where $\dot{\BPhi}_i$ is a tangent vector of $\Od(d)$ at $\BPhi_i$.
All the derivations can be found in Lemma~\ref{label: gradient and hessian}.
Thus, the Newton step at $\BZ_{i}$ is given by
\begin{equation} \label{eq:newtonstep}
		\widetilde{\BV}_{i} :=  \BH_{i} ^{-1}\lp \sum_{j=1}^n \BC_{ij}\widehat{\BZ}_{j} - \BH_{i}\BZ_{i}\rp = \frac{1}{2}\BH_i^{-1} \sum_{j=1}^n \left(\BC_{ij}\widehat{\BZ}_j\BZ_i^{\top} - \BZ_i \widehat{\BZ}_j^{\top}\BC_{ji} \right) \BZ_i 
		\end{equation}
where 
\begin{equation}\label{def:HI}
\BH_{i} := \frac{1}{2}\sum_{j=1}^n \lp \BC_{ij}\widehat{\BZ}_{j}\BZ_{i}^{\top} + \BZ_{i}\widehat{\BZ}^{\top}_{j}\BC_{ji}\rp
\end{equation}
which is associated with the Hessian of $f(\cdot)$ at $\BZ_i.$
In fact, we show in Proposition~\ref{prop:2nderror} that $\widetilde{\BV}_i$ minimizes the quadratic approximation of $f(\cdot)$ at $\BZ_i$ on the tangent space. The  one-step Newton approximation of the MLE is then given by $\widetilde{\BZ}_i := \PP(\BZ_i + \widetilde{\BV}_i)$.

Note that the distribution of $\BZ_i + \widetilde{\BV}_i$ is relatively simpler to characterize, and then it remains to upper bound the approximation error $\|\BE_i\|_F:$ 
\begin{equation}\label{eq:error_decomp}
 \BE_i : =  \widehat{\BZ}_i - (\BZ_i + \widetilde{\BV}_i )\end{equation}
which can be decomposed as
\begin{equation} \label{eqn: error decomposition}
\begin{aligned}
	\norm{\BE_{i}}_{F} 
	&= \norm{\widehat{\BZ}_{i} - \widetilde{\BZ}_i + \widetilde{\BZ}_i- (\BZ_{i} + \widetilde{\BV}_{i}) }_{F} \le \norm{\widehat{\BZ}_{i} - \widetilde{\BZ}_{i}}_{F} + \|\BZ_i + \widetilde{\BV}_i - \PP(\BZ_i + \widetilde{\BV}_i)\|_F.
\end{aligned}
\end{equation}

The challenging part of estimating~\eqref{eqn: error decomposition} arises from $\|\widehat{\BZ}_i - \widetilde{\BZ}_i\|_F$. We employ the following approach: let $\widehat{\BV}_i$ be the tangent vector at $\BZ_i$ such that $\widehat{\BZ}_i = \PP(\BZ_i + \widehat{\BV}_i)$. The existence of such $\widetilde{\BV}_i$ is guaranteed if $\|\widehat{\BZ}_i - \BZ_i\| < 1$ by Lemma~\ref{lemma: Lipschitz}. Then, since the polar retraction $\PP(\cdot)$ is essentially a projection onto the convex set $\{\BY: \|\BY\|\leq 1\}$ for any input with operator norm greater than 1, we can verify that $\mathcal{P}(\cdot)$ is $1$-Lipschitz continuous: 
\[
\|\widehat{\BZ}_i - \widetilde{\BZ}_i\|_F = \|\PP(\BZ_{i} + \widehat{\BV}_i) - \PP(\BZ_{i} +\widetilde{\BV}_i)\|_F \leq \|\widehat{\BV}_i - \widetilde{\BV}_i\|_F.
\]
As a result, it suffices to estimate $\|\widehat{\BV}_i - \widetilde{\BV}_i\|_F$ and $\|\BZ_i + \widetilde{\BV}_i - \PP(\BZ_i + \widetilde{\BV}_i)\|_F$, which are given by the following results.

\begin{proposition}[\bf Distribution of $\widetilde{\BV}_i$]\label{prop:tildeV}
Under the conditions of Theorem~\ref{thm:main}, for every $1\leq i\leq n$, 
\begin{align*}
& \widetilde{\BV}_i = \frac{\sigma \sqrt{n-1}}{2n} (\BG_i - \BG_i^{\top}) \BZ_i+ O \lp\SNR^{-2}(\sqrt{d}+\sqrt{\log n})^2 \rp, \\
& \|\widetilde{\BV}_i\|_F \lesssim \SNR^{-1} (\sqrt{d} + \sqrt{\log n}),
\end{align*}
where $\BG_{i}$ is a $d\times d$ Gaussian random matrix with i.i.d. $\mathcal{N}(0,1)$ entries.
\end{proposition}

\begin{proposition}[\bf Distance between MLE and one-step Newton iteration]\label{prop:newtonstep}
Under the conditions of Theorem~\ref{thm:main}, the following statements hold with probability at least $1-n^{-2}$:
\begin{align*}
& \|\widehat{\BV}_i\| \lesssim \SNR^{-1} (\sqrt{d} + \sqrt{\log n}), \\
& \norm{\widehat{\BZ}_{i} - \widetilde{\BZ}_{i}}_{F} \leq \|\widehat{\BV}_i - \widetilde{\BV}_i\|_F \lesssim \SNR^{-2}(\sqrt{d} + \sqrt{\log n})^2.
\end{align*}
\end{proposition}

\begin{lemma}[{\cite[Appendix E.1]{LSW19}}]\label{lem:retract}
Suppose $\| \BV_i \|_F\leq 1$, and then 
\[
 \|\BZ_i + \BV_i - \PP(\BZ_i + \BV_i)\|_F \leq \|\BV_i\|_F^2.
\]
\end{lemma}

With Proposition~\ref{prop:tildeV} and~\ref{prop:newtonstep}, we are ready to give a proof of Theorem~\ref{thm:main}. The proof of Proposition~\ref{prop:tildeV} and~\ref{prop:newtonstep} can be found in Section~\ref{ss:propkey}.

\begin{proof}[\bf Proof of Theorem~\ref{thm:main}] 
Applying Proposition~\ref{prop:newtonstep} and Lemma~\ref{lem:retract} gives:
\[
 \|\BZ_i + \widetilde{\BV}_i - \PP(\BZ_i + \widetilde{\BV}_i)\|_F  \leq \|\widetilde{\BV}_i\|_F^2 \lesssim \SNR^{-2}(\sqrt{d} + \sqrt{\log n})^2
\]
and thus
\begin{align*}
\|\BE_i\|_F & \leq \|\widehat{\BV}_i - \widetilde{\BV}_i\|_F + \|\widetilde{\BV}_i\|_F^2  \lesssim  \SNR^{-2}(\sqrt{d} + \sqrt{\log n})^2 
\end{align*}
provided that $\SNR\gtrsim \sqrt{d} + \sqrt{\log n}.$ Note that Proposition~\ref{prop:tildeV} implies
\[
\widetilde{\BV}_i = \frac{\sigma \sqrt{n-1}}{2n} (\BG_i -\BG_i^{\top}) \BZ_i+ O \lp \SNR^{-2}(\sqrt{d} + \sqrt{\log n})^2 \rp
\]
and therefore, it holds
\begin{align*}
\widehat{\BZ}_i & = \BZ_i + \widetilde{\BV}_i + \BE_i  = \BZ_i + \frac{\sigma \sqrt{n-1}}{2n} (\BG_i - \BG_i^{\top}) \BZ_i+ O\left( \SNR^{-2}(\sqrt{d} + \sqrt{\log n})^2 \right)
\end{align*}
which finishes the proof.
\end{proof}

\subsection{Proof sketch of Theorem~\ref{thm:spectral_main}: distribution of spectral estimator}
Similar to~\eqref{eqn: local MLE}, we will argue that
\begin{equation}\label{eq:local spectral}
\widehat{\BS}_i^{\eig} = \argmin_{\BS_i^{\top}\BS_i = \BPi_i } -\sum_{j=1}^n \lag \BC_{ij}\widehat{\BS}_j^{\eig}, \BS_i\rag
\end{equation}
where $\widehat{\BS}^{\rm eig}$ is the global minimizer to~\eqref{eq:spectral} and $\BPi_i := \widehat{\BS}^{\rm eig\top}_{i}\widehat{\BS}^{\rm eig}_i$ satisfies $\sum_{i=1}^n \BPi_i = n\I_d$. 
The proof follows from Lemma~\ref{eq:spr_local}.

Since our goal is to derive the distribution of $\PP(\BSeig_i)$, it will be more reasonable to directly analyze the rounded $\BSeig_i$. This motivates our transformation of~\eqref{eq:local spectral} to a scenario that is similar to~\eqref{eqn: local MLE} where the constraint is $\Od(d)$. Let $\BPhi = \BS_i \BPi_i^{-1/2}$ as $\BPi_i$ is full-rank. Then minimizing~\eqref{eq:local spectral} is equivalent to the following optimization:
\begin{equation} \label{eq:spectral equivalent}
\min_{\BS_i^{\top}\BS_i = \BPi_i } -\sum_{j=1}^n \lag \BC_{ij}\widehat{\BS}_j^{\eig}, \BS_i\rag = 	\underset{\quad \BPhi \in \Od(d)}{\min} \quad - \sum_{j=1}^n\lag \BC_{ij}\widehat{\BS}^{\rm eig}_j\BPi_i^{1/2}, \BPhi \rag
\end{equation}
and Lemma~\ref{eq:spr_local} implies
\[
\widehat{\BZ}_i^{\eig} : = \PP(\widehat{\BS}_i^{\eig}) = \argmin_{ \BPhi\in \Od(d)} \quad - \sum_{j=1}^n\lag \BC_{ij}\widehat{\BS}^{\rm eig}_j\BPi_i^{1/2}, \BPhi\rag.
\]

This resembles the formulation in the scenario of MLE.
Given the transformation in~\eqref{eq:spectral equivalent}, the error term $\BE^{\rm eig}_{i} :=  \BZeig - (\BZ_i + \widetilde{\BV}^{\rm eig}_i )$ will depend on the conditioning of $\BPi_{i}$. We define the condition number as $\kappa := \pi_{\max}/\pi_{\min} \ge 1$
where  Lemma~\ref{lem:pi} estimates the largest and smallest eigenvalue of $\BPi_{i}$ via
$\|\BPi_i - \I_d\| \lesssim \SNR^{-1} d^{-1/2}(\sqrt{d} + \sqrt{\log n}).$ 
Then we are ready to present the following results similar to those in Section~\ref{ss:proof sketch mle}.
\begin{proposition}[\bf Distribution of $\widetilde{\BV}^{\rm eig}_i$]\label{prop:tildeV_spectral}
Under the conditions of Theorem~\ref{thm:spectral_main}, for every $1\leq i\leq n$, 
\begin{align*}
& \widetilde{\BV}^{\rm eig}_i = \frac{\sigma \sqrt{n-1}}{2n} (\BG_i - \BG_i^{\top}) \BZ_i+ O \lp \sqrt{\kappa} \SNR^{-2}(\sqrt{d}+\sqrt{\log n})^2 \rp, \\
& \|\widetilde{\BV}_i^{\eig}\|_F\lesssim \SNR^{-1} (\sqrt{d} +\sqrt{\log n})
\end{align*}
where $\BG_{i}$ is a $d\times d$ Gaussian random matrix with i.i.d. $\mathcal{N}(0,1)$ entries.
\end{proposition}

\begin{proposition}[\bf Distance between spectral estimator and one-step Newton iteration]\label{prop:newtonstep_spectral}
Under the conditions of Theorem~\ref{thm:spectral_main}, the following statements hold with high probability:
\begin{align*}
& \|\widehat{\BV}_i^{\rm eig}\|  \lesssim \SNR^{-1}(\sqrt{d}+\sqrt{\log n}), \\
& \norm{\BZeig_{i} - \widetilde{\BZ}^{\eig}_{i}}_{F}  \leq \|\widehat{\BV}_i^{\rm eig} - \widetilde{\BV}^{\rm eig}_i\|_F \lesssim \kappa^{1/2}\SNR^{-2}(\sqrt{d}+\sqrt{\log n})^2,
\end{align*}
where $\widetilde{\BZ}_i^{\eig} = \PP(\BZ_i + \widetilde{\BV}^{\rm eig}_i)$. 
\end{proposition}

The proof of Theorem~\ref{thm:spectral_main} follows from
the exact same argument as the Proof of Theorem~\ref{thm:main}: combining Proposition~\ref{prop:newtonstep_spectral} with Lemma~\ref{lem:retract} gives an upper bound on 
\begin{align*}
\norm{\BZeig_i - (\BZ_i + \widetilde{\BV}^{\rm eig}_i )}_{F} 
& \leq \left\| \BZeig_i - \widetilde{\BZ}^{\eig}_{i} \right\|_F+ \left\|\PP(\BZ_i + \widetilde{\BV}^{\rm eig}_i) - (\BZ_i + \widetilde{\BV}^{\rm eig}_i ) \right\|_F \\
& \leq \| \widehat{\BV}_i^{\eig} - \widetilde{\BV}_i^{\eig} \|_F + \|\widetilde{\BV}_i^{\eig}\|_F^2 \lesssim \sqrt{\kappa}\SNR^{-2}(\sqrt{d} +\sqrt{\log n})^2,
\end{align*}
and the distribution of $\BZ_i + \widetilde{\BV}_i^{\eig}$ follows from Proposition~\ref{prop:tildeV_spectral}. 

\section{Numerics}\label{s:numerics}
In this section, we will present the experimental results that empirically confirm Theorem~\ref{thm:main} and~\ref{thm:spectral_main}. 
We select $n = 800$ and $d = 5$, and generate data at various noise levels $\sigma = C_{\sigma} \sqrt{n/d}$, with $C_{\sigma}$ ranging from 0.05 to 0.5 in increments of 0.05. By definition, the signal-to-noise ratio ($\SNR$) satisfies 
$C_{\sigma} = \SNR^{-1}$ and thus  $\SNR$ is between $2$ and $20$.
 For each noise level, we compute the spectral estimator $\BZeig = \PP_{n}(\BSeig)$ where $\BSeig$ is the top-$d$ eigenvectors of $\BC$ subject to $\BSeigt\BSeig = n \BI_{d}$.   Since solving~\eqref{eq:sdp} is computationally expensive, we take an alternative way to find the SDP solution by applying the generalized power method: $\widehat{\BZ}^{t+1} = \PP_{n}(\BC\widehat{\BZ}^{t})$ for $t=0,1,2,\cdots, L-1$ with $L=100$ and spectral initialization $\widehat{\BZ}^{0} = \BZeig$. We set $\widehat{\BZ} = \widehat{\BZ}^L$ if
$\| (\BLambda - \BC) \widehat{\BZ} \|_F \leq 10^{-8}, \lambda_{d+1}(\BLambda - \BC) \geq 10^{-8},$
 i.e., $\widehat{\BZ}$ satisfies the approximate KKT condition of the~\eqref{eq:sdp}.

There are a few theoretical aspects we will verify in numerical simulations. Note that Theorem~\ref{thm:main} and~\ref{thm:spectral_main} imply that
\begin{equation}\label{eq:errorterm}
\widehat{\BZ}_{i}\BQ^{\top}\BZ_i^{\top}  = \I_d + \frac{\sigma\sqrt{n-1}}{\sqrt{2}n} \BG_i^{\sk} + O \lp \SNR^{-2}(\sqrt{d} + \sqrt{\log n})^2 \rp,
\end{equation}
where  $\BQ = \PP(\BZ^{\top}\widehat{\BZ})$ and $\BG_i^{\sk} := (\BG_i - \BG_i^{\top})/\sqrt{2}$
is a skew-symmetric standard Gaussian random matrix, and a similar result holds for spectral estimator. We will test (a) whether the result above is indeed a second-order expansion of the MLE and spectral estimator; (b) the first-order error is nearly Gaussian; (c) the coefficient in front of $\BG_i^{\sk}$ is $\sigma\sqrt{n-1}/(\sqrt{2}n).$

 To verify the main results, this section includes boxplots of test statistics for each noise level, with each box aggregating results from $n=800$ blocks. The actual test statistic depends on the settings. The boxplot is structured conventionally:  the central line represents the median, the edges of the box represent the $25$-th and $75$-th percentiles, the whiskers extend to the smallest and largest values within $1.5$ times the interquartile range (IQR) from the quartiles, and individual points outside the whiskers are outliers. 

\subsection{Verification of the second-order error}

\begin{figure}[H] 
    \centering
    \begin{minipage}{0.45\textwidth}
        \centering
        \includegraphics[width=\linewidth]{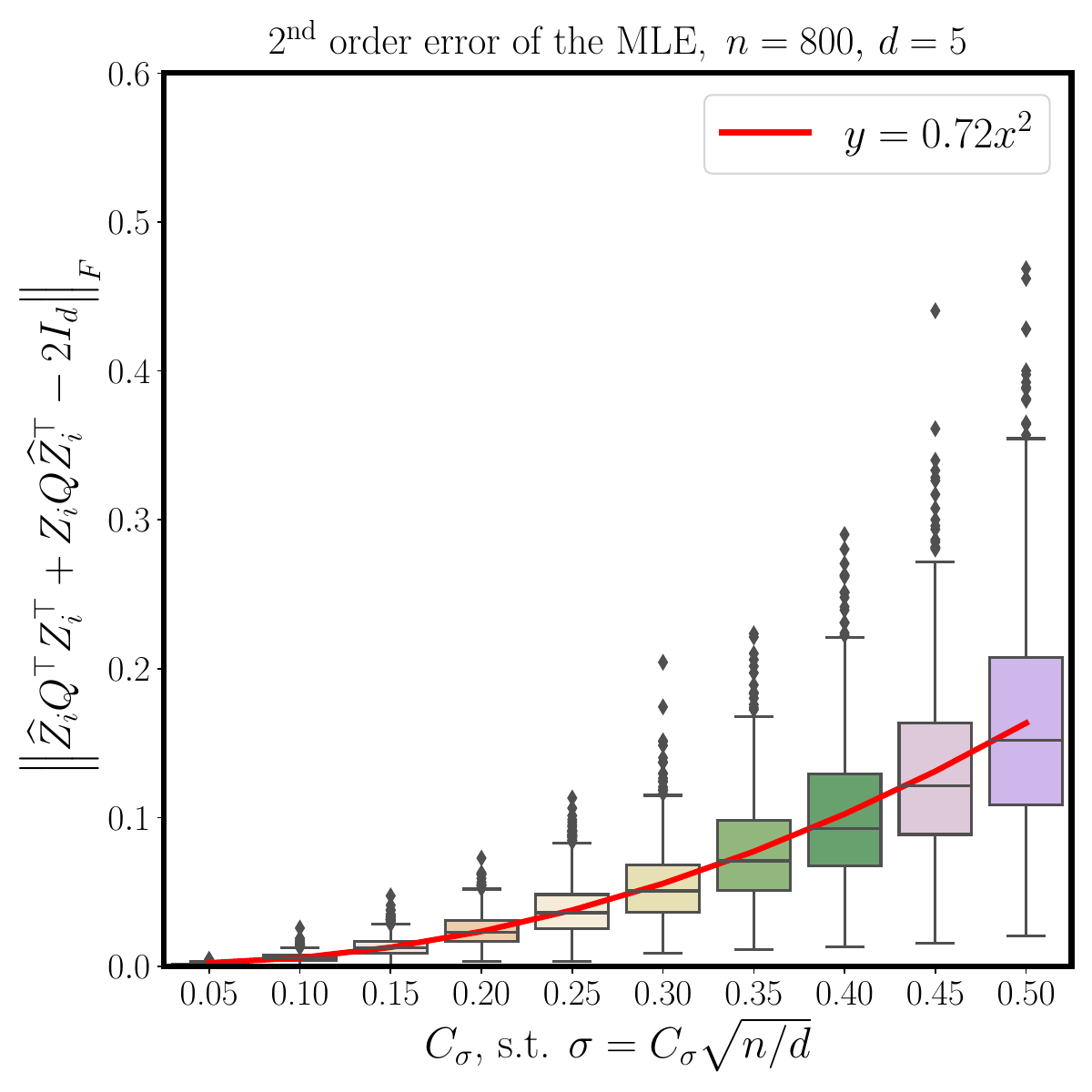} 
    \end{minipage}
    \hfill
    \begin{minipage}{0.45\textwidth}
        \centering        \includegraphics[width=\linewidth]{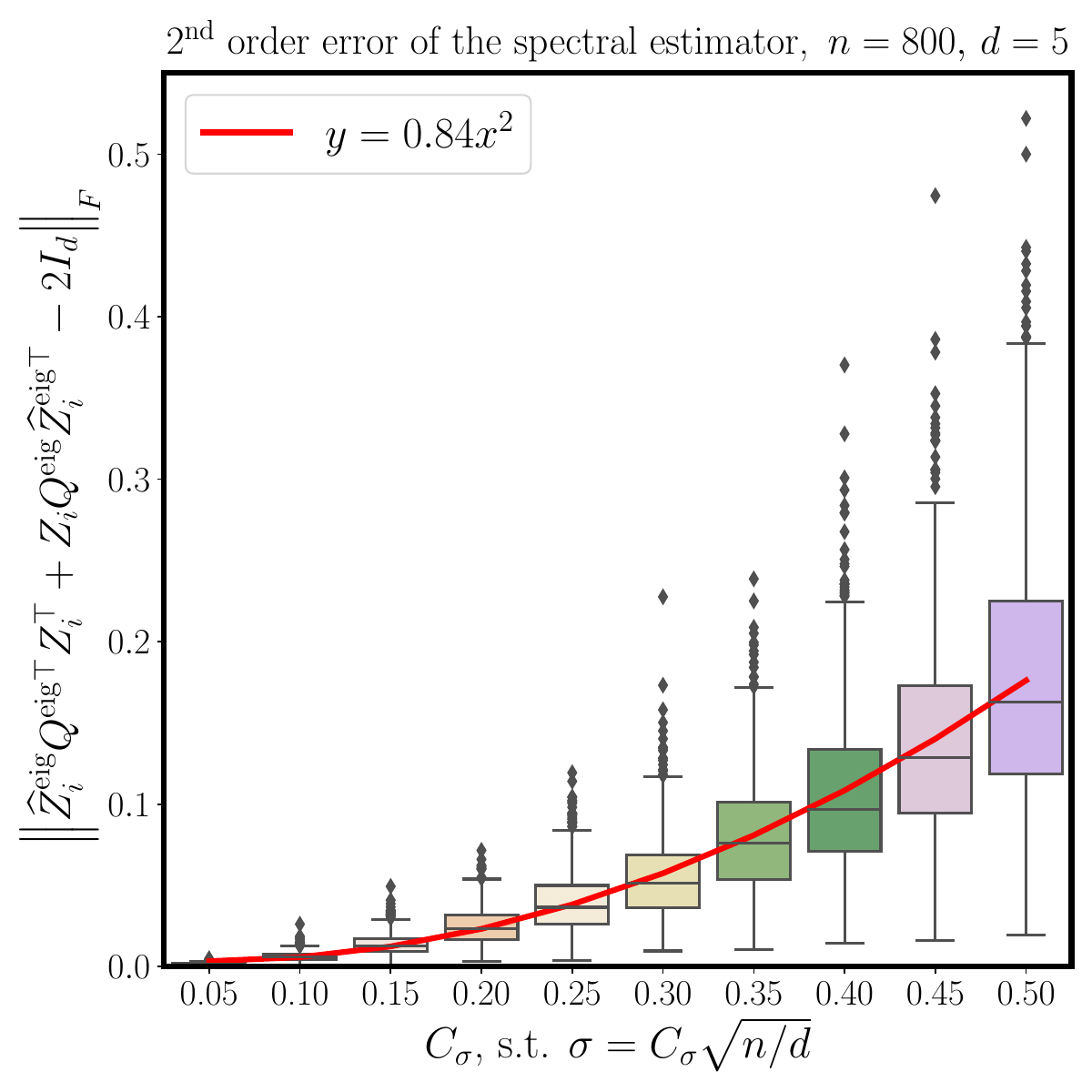} 
    \end{minipage}
    \caption{Left: MLE; right: spectral estimator. The figure confirms the second-order error term in the expansion of~\eqref{eq:main} and~\eqref{eq:main_spectral}.}
    \label{fig:quadratic}
\end{figure}

Based on Theorem~\ref{thm:main} and~\ref{thm:spectral_main}, we have
$\| \widehat{\BZ}_{i}\BQ^{\top}\BZ_i^{\top} + \BZ_i\BQ \widehat{\BZ}_i^{\top} -2\I_d\|_F = O(\SNR^{-2}(\sqrt{d} + \sqrt{\log n})^2)$
where $\SNR^{-1} = C_{\sigma}$ and the skew-symmetric term is cancelled.
As shown in Figure~\ref{fig:quadratic}, when $\sigma$ increases, the error  scales quadratically w.r.t. $\SNR^{-1}$ for both estimators. It also implies that the first-order error term, i.e.,  a skew-symmetric term, is cancelled, which confirms our theory. 
Moreover, we observe that the second-order error of the spectral estimator is slightly larger than that of MLE, as suggested by the coefficients of the fitted curves ($0.84>0.72$). Even though the first-order error is the same for both estimators, the MLE has a smaller second-order error which may serve as an evidence for the existence of $\kappa > 1$ in~\eqref{eq:main_spectral}.

\subsection{Test of the approximate normality}
\begin{figure}[h!]
    \centering
    \begin{minipage}{0.45\textwidth}
        \centering
        \includegraphics[width=\linewidth]{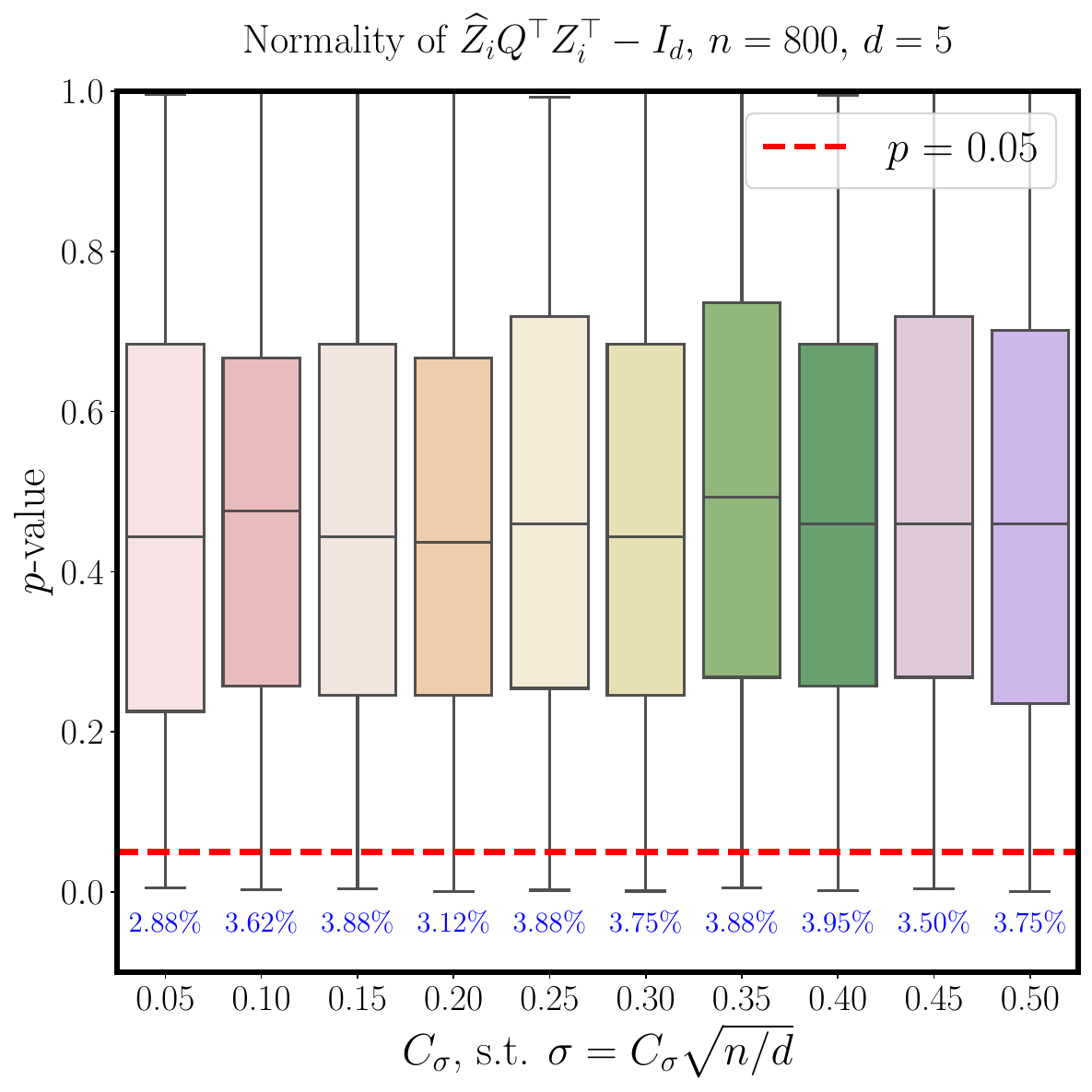} 
    \end{minipage}
    \hfill
    \begin{minipage}{0.45\textwidth}
        \centering
        \includegraphics[width=\linewidth]{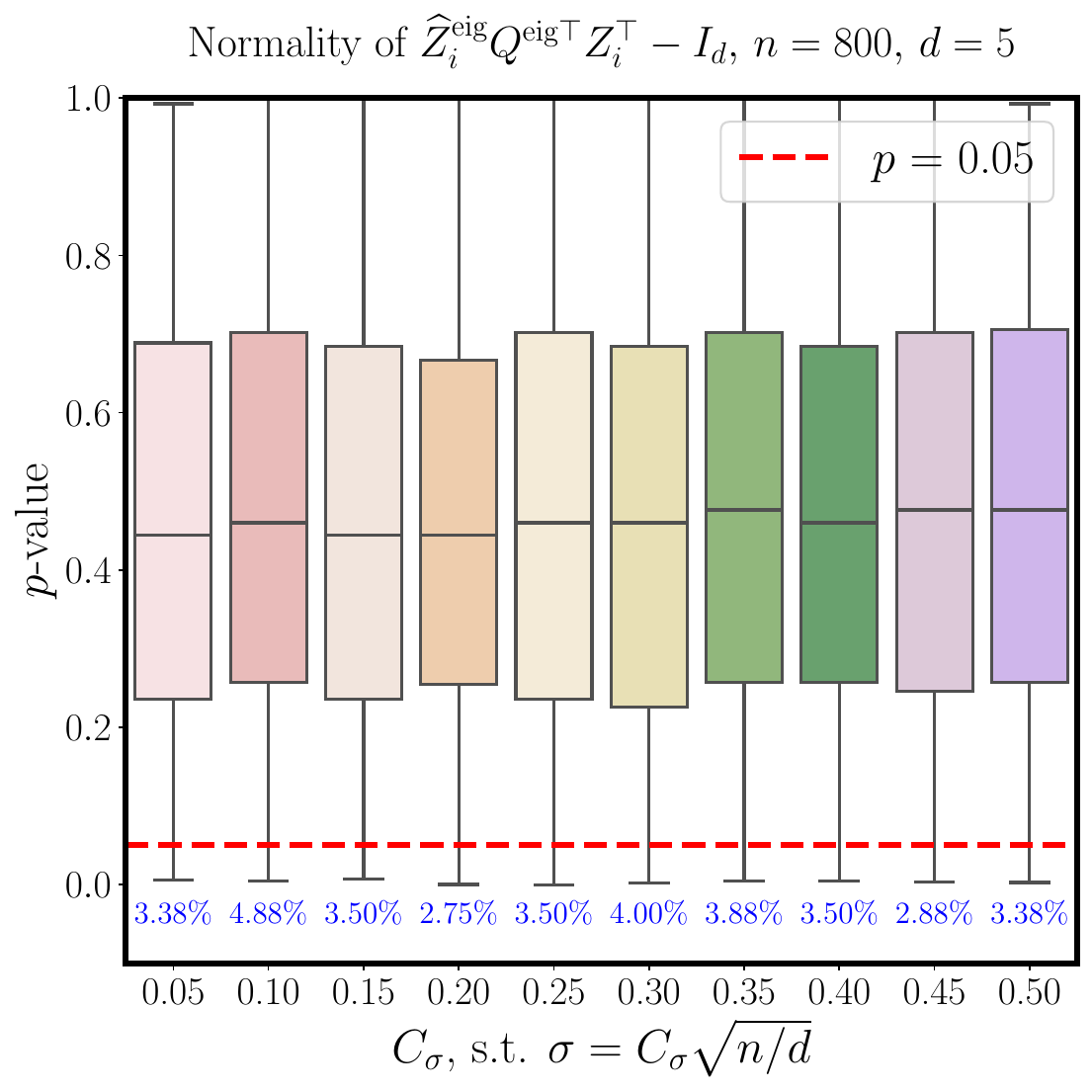} 
    \end{minipage}
    \caption{Left: MLE; right: spectral estimator. Here is the boxplot of $p$-values from the Kolmogorov-Smirnov test assessing the normality of the dominant terms for both MLE and spectral estimators across various noise levels.}
    \label{fig:Gaussianality}
\end{figure}

We proceed to test the normality of the first-order error $\widehat{\BZ}_i (\BZ_i\BQ)^{\top} -\I_d$. Our theorems imply that the upper triangular entries of $\frac{\sqrt{2}n}{\sigma \sqrt{n-1}} \cdot (\widehat{\BZ}_i (\BZ_i\BQ)^{\top} -\I_d)$ are i.i.d. standard Gaussian variables. To test the entries' normality, we compute  the Kolmogorov-Smirnov statistic for each $1\leq i\leq 800$ and also calculate the $p$-value under the null hypothesis that the upper triangle entries are i.i.d. standard normal.
We compile the $p$-values from $n = 800$ blocks into a single box in the left subplot of Figure~\ref{fig:Gaussianality} for each noise level and carry out the same procedure for the spectral estimator, with results displayed in the right subplot of Figure~\ref{fig:Gaussianality}. The proportions of $p$-values below $0.05$ in each box are reported below each error bar. The figure indicates that there is no strong evidence against the  normality of the leading error term in Theorem~\ref{thm:main} and~\ref{thm:spectral_main}.


\subsection{Determine the coefficient of the skew-symmetric Gaussian terms}
\begin{figure}[h!]
    \centering
    \begin{minipage}{0.45\textwidth}
        \centering
        \includegraphics[width=\linewidth]{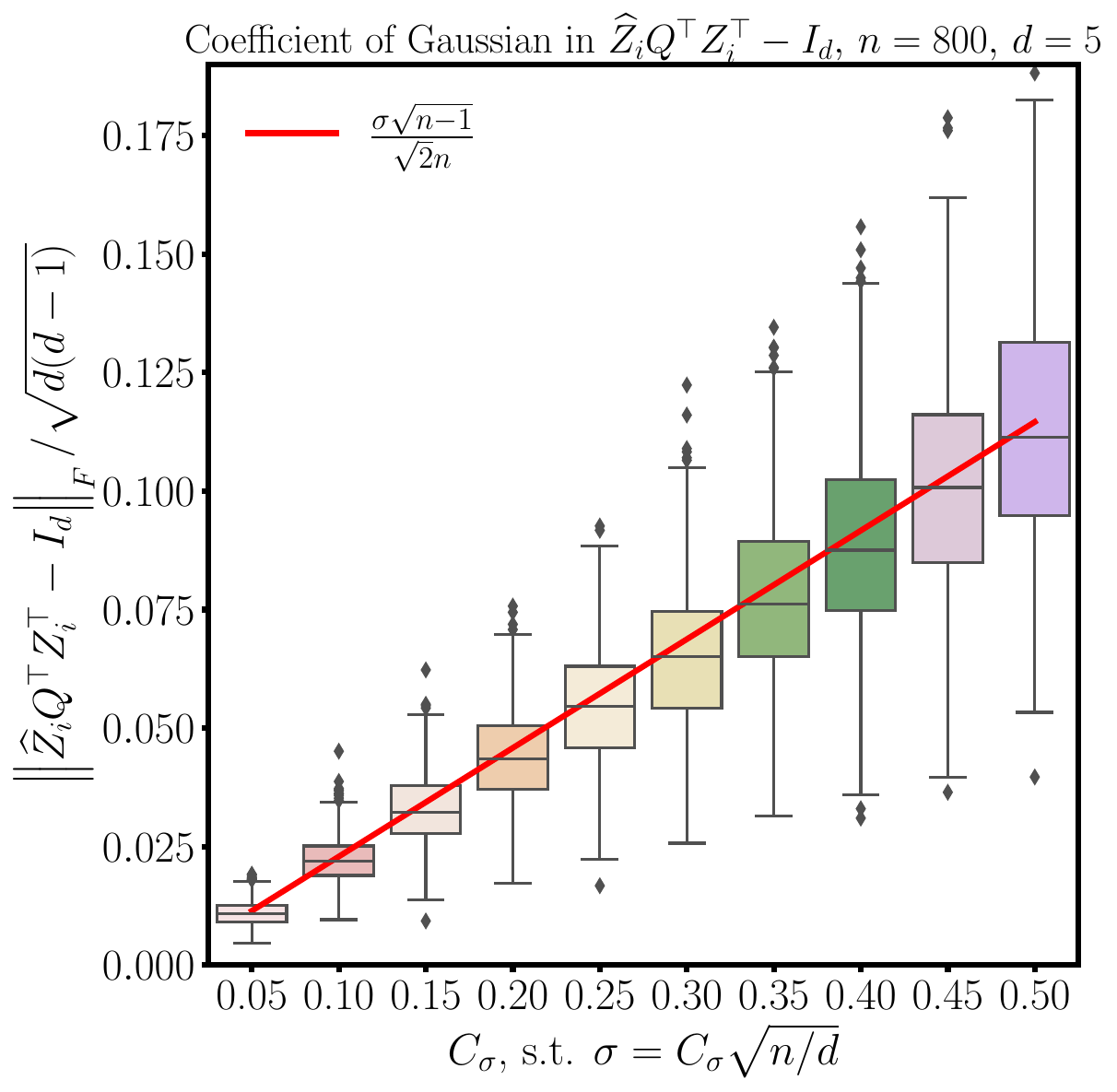} 
    \end{minipage}
    \hfill
    \begin{minipage}{0.45\textwidth}
        \centering
        \includegraphics[width=\linewidth]{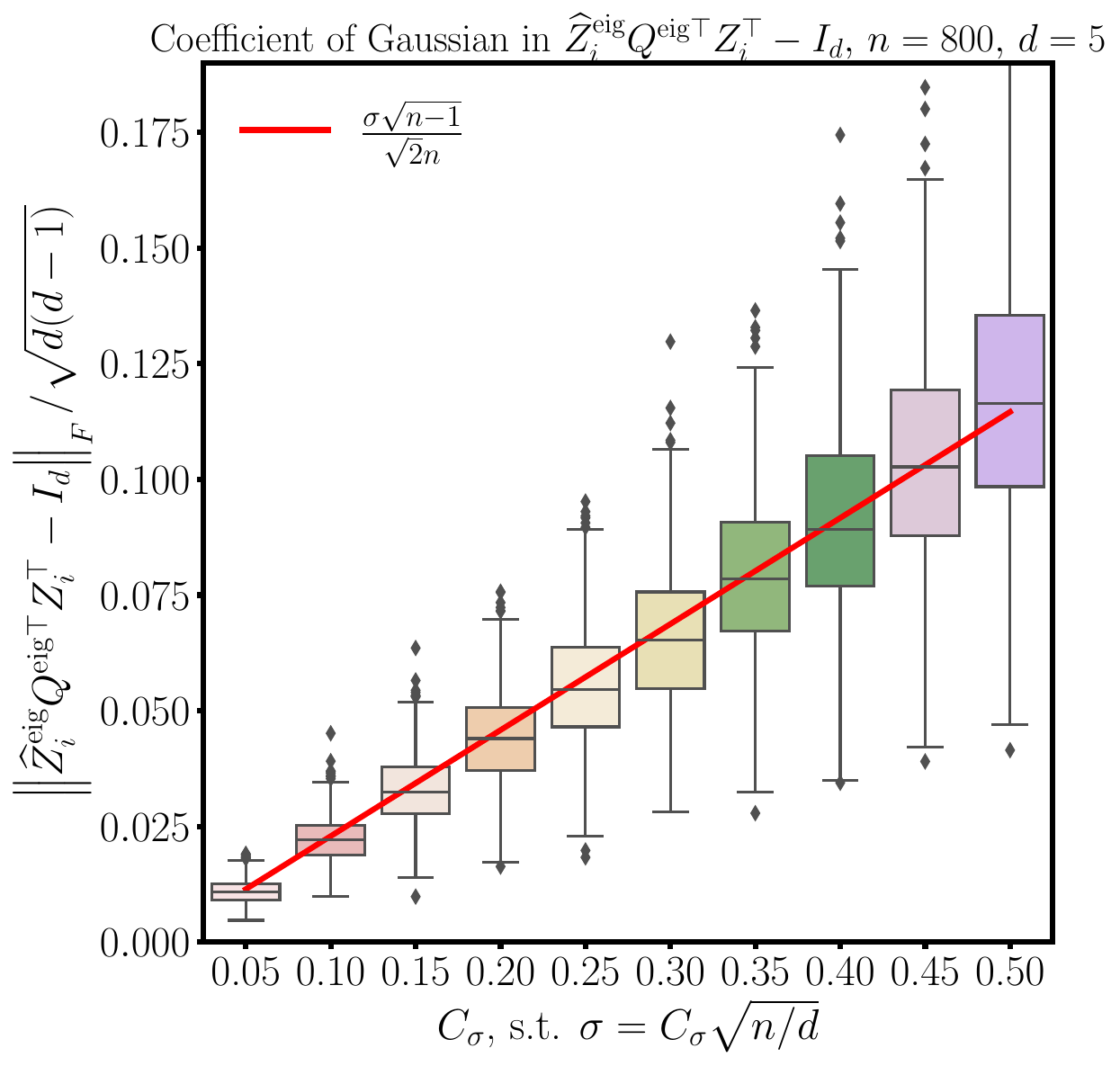} 
    \end{minipage}
    \caption{Left: MLE; right: spectral estimator. The theoretical value of the coefficient matches the experimental data as $\sigma$ increases, verifying the coefficients of the Gaussian term in~\eqref{eq:main} and~\eqref{eq:main_spectral}.}
    \label{fig:coef of gaussian}
\end{figure}

Finally, we verify the coefficient of the skew-symmetric Gaussian term $\BG_i^{\sk}$ in~\eqref{eq:errorterm} is equal to $\sigma \sqrt{n-1}/(\sqrt{2}n)$. As we do not know $\|\BG_i^{\sk}\|_F^2$, we approximate it by its expectation $d(d-1)$. 
For both MLE and the spectral estimator under each noise level, we calculate 
$\| \HBZ_{i} - \BZ_{i}\BQ\|_{F}/\sqrt{d(d-1)}$
and plot the boxplot for each $\sigma$. Then we compare it with the theoretical value $\sigma \sqrt{n-1} / (\sqrt{2}n)$. As shown in Figure~\ref{fig:coef of gaussian}, the theoretical value nearly matches the median of each error bar as $\sigma$ increases. This is a strong evidence to show our theoretical prediction is aligned with the coefficients of the skew-symmetric Gaussian term in~\eqref{eq:errorterm}, and thus~\eqref{eq:main} and~\eqref{eq:main_spectral}.

\appendix

\section{Riemannian gradient and Hessian}\label{s:riemann}
We will derive the Riemannian gradient and Hessian of functions of the form
\begin{equation}\label{eq:f_stiefel}
f(\BPhi) = \left \lag \BA, \BPhi \right \rag
\end{equation}
where $\BA \in \reals^{d \times p}$ is arbitrary and
\[
	\BPhi \in \St(d,p)=\{\BPhi \in \reals^{p\times d}: \BPhi\BPhi^{\top} = \BI_{d},~d\le p \}
\]
which is the Stiefel manifold.
Consider the optimization
\[
\min_{\BPhi\in\St(d,p)}~f(\BPhi)
\]
The global minimizer is given by $\widehat{\BPhi} = -(\BA\BA^{\top})^{-1/2}\BA$ if $\BA$ is full row rank and
$f(\widehat{\BPhi}) = - \|\BA\|_*$ 
where $\|\BA\|_*$ is the nuclear norm of $\BA.$
The following lemma gives the Riemannian gradient and Hessian of~\eqref{eq:f_stiefel} at $\BPhi \in \St(d,p)$.

\begin{lemma} \label{label: gradient and hessian}
The Riemannian gradient and the quadratic form associated to the Riemannian Hessian of~\eqref{eq:f_stiefel} at $\BPhi \in \St(d,p)$ are given by
\begin{align*}
\rgrad f(\BPhi) & = \BA - \frac{\BA\BPhi^{\top} + \BPhi\BA^{\top}}{2}\BPhi, \\
\Hess f(\BPhi) [\dot{\BPhi},\dot{\BPhi}] & = -\frac{1}{2} \lag\BA\BPhi^{\top} +\BPhi\BA^{\top}, \dot{\BPhi}\dot{\BPhi}^{\top} \rag,
\end{align*}
where $\dot{\BPhi}$ is an element in $T_{\BPhi}(\St(d,p))$.
\end{lemma}
\begin{proof}

Let $\Proj_{\BPhi}(\cdot)$ be the projection operator from the Euclidean space to the tangent space of $\St(d,p)$ at $\BPhi$, i.e., 
\[
T_{\BPhi}(\St(d,p)) : = \left \{ \dot{\BPhi} \in \reals^{d \times p} :\BPhi\dot{\BPhi}^{\top} + \dot{\BPhi}\BPhi^{\top} = 0 \right \}
\] 
and the projection is
\[
	\Proj_{\BPhi}(\BY) : = \BY - \frac{1}{2}(\BPhi\BY^{\top} + \BY\BPhi^{\top})\BPhi
\] 
for any $\BY\in\RR^{d\times p}.$
The Riemannian gradient is given by
\begin{equation} 
	\begin{aligned}
	\rgrad f(\BPhi ) = \Proj_{\BPhi}\lp\grad f(\BPhi )\rp = \Proj_{\BPhi}\lp \BA\rp = \BA - \frac{\BA\BPhi^{\top} + \BPhi\BA^{\top}}{2}\BPhi,
\end{aligned} 
\end{equation}
where $\nabla f(\BPhi) = \BA$ is the gradient of $f(\cdot)$ in the ambient space. We proceed to compute the Hessian of $f(\BPhi).$
The directional gradient of $\rgrad f(\BPhi )$ in the direction of $\dot{\BPhi}$ at $\BPhi$ is given by
\begin{align*}
\text{D} \rgrad f(\BPhi )[\dot{\BPhi}] & : = \lim_{t\rightarrow 0} \frac{1}{t} \left[\BA -  \frac{\BA(\BPhi+t\dot{\BPhi})^{\top} + (\BPhi+t\dot{\BPhi})\BA^{\top}}{2}(\BPhi+t\dot{\BPhi}) - \left( \BA -  \frac{\BA\BPhi^{\top} + \BPhi\BA^{\top}}{2}\BPhi\right) \right] \\
& = \frac{1}{2} \lim_{t\rightarrow 0} \frac{1}{t} \left[ - (\BA(\BPhi+t\dot{\BPhi})^{\top} + (\BPhi+t\dot{\BPhi})\BA^{\top})(\BPhi+t\dot{\BPhi}) + \left( \BA\BPhi^{\top} + \BPhi\BA^{\top}\right)\BPhi \right] \\
& = \frac{1}{2} \lim_{t\rightarrow 0} \frac{1}{t} \left[ - (\BA\BPhi^{\top} +\BPhi\BA^{\top} + t(\BA\dot{\BPhi}^{\top} + \dot{\BPhi}\BA^{\top}))(\BPhi+t\dot{\BPhi}) + \left( \BA\BPhi^{\top} + \BPhi\BA^{\top}\right)\BPhi \right] \\
& = - \frac{\BA\BPhi^{\top} +\BPhi\BA^{\top}}{2}\dot{\BPhi} - \frac{1}{2}(\BA\dot{\BPhi}^{\top} + \dot{\BPhi}\BA^{\top}) \BPhi.
\end{align*}
Then the quadratic form associated to $\Hess f(\BPhi)$ is given by
\[
\Hess f(\BPhi) [\dot{\BPhi},\dot{\BPhi}]= \left \lag \text{D} \rgrad f(\BPhi )[\dot{\BPhi}], \dot{\BPhi}\right \rag = -\frac{1}{2} \lag\BA\BPhi^{\top} +\BPhi\BA^{\top}, \dot{\BPhi}\dot{\BPhi}^{\top} \rag
\]
where $\dot{\BPhi}$ is an element on the tangent space satisfying $\BPhi\dot{\BPhi}^{\top} + \dot{\BPhi}\BPhi^{\top} = 0.$
\end{proof}

In our analysis, we need to compute the Newton step at a given initial value $\BPhi_{0}.$ Now we will investigate the property of the one-step Newton iteration of $f(\BPhi)$ at $\BPhi_{0}\in \St(d,p)$.

\begin{definition}[Polar retraction]
For any full row rank matrix $\BY \in \reals^{d \times p}$, its projection onto the Stiefel manifold is given by
\[
\mathcal{P}(\BY) := (\BY\BY^{\top})^{-1/2}\BY.
\] 
For any $\BPhi_{0}\in \St(d,p)$ and any $\BV$ on the tangent space of $\St(d,p)$ at $\BPhi_0$, i.e., $\BV\BPhi_{0}^{\top} +\BPhi_{0}\BV^{\top}= 0$, then the corresponding retraction of $\BPhi_0 + \BV$ is given by
\begin{align}\label{eqn: polar retraction}
R_{\BPhi_0}(\BV) := \mathcal{P}(\BPhi_0 + \BV) = [(\BPhi_{0} + \BV)(\BPhi_{0} + \BV)^{\top}]^{-\frac{1}{2}}(\BPhi_{0} + \BV)  = (\I_d + \BV\BV^{\top})^{-\frac{1}{2}}(\BPhi_{0}+\BV)
\end{align}
where $\BPhi_0 + \BV$ is full row-rank and $\BPhi_0\BPhi_0^{\top} = \I_d.$
\end{definition}

Next, we will analyze the distance between the one-step Newton iteration and the global minimizer. Let $\widehat{\BPhi}$ be the global minimizer to $f$, and let $\widehat{\BV}$ be the corresponding tangent vector such that
\[
R_{\BPhi_0}(\widehat{\BV}) = \widehat{\BPhi} ~\Longleftrightarrow~
(\I_d + \widehat{\BV}\widehat{\BV}^{\top})^{-1/2}(\BPhi_{0}+\widehat{\BV}) = \widehat{\BPhi}.
\]
In particular, if $\widehat{\BPhi}$ is close to $\BPhi_{0}$, i.e., $\|\widehat{\BPhi} - \BPhi_0\| < 1$, then such a $\widehat{\BV}$ always exists. Moreover, the following lemmas will characterize the Lipschitz continuity of $R_{\BPhi_0}(\cdot)$ and its invertibility. 
\begin{lemma}[\cite{LSW19,LCD+19}] \label{lem:retraction}
	Let $R_{\BPhi_0}(\cdot)$ be the polar retraction defined in~\eqref{eqn: polar retraction} and $\norm{\BV}_{F} \le 1$, then for any $\BY \in \St(d,p)$, we have 
	\[
	\norm{R_{\BPhi_{0}}(\BV)-\BY}_{F} \le \norm{\BPhi_{0}+\BV- \BY}_{F}.
	\]
\end{lemma}
We omit the proof here as it has been introduced in~\cite{LSW19,LCD+19}. The idea of why Lemma~\ref{lem:retraction} holds follows from the fact that ${\cal P}(\cdot)$ is essentially the projection  onto the unit ball of operator norm for any matrix with $\|\cdot\| > 1.$ Then
\[
\norm{R_{\BPhi_{0}}(\BV)-\BY}_{F} = \norm{\PP(\BPhi_0 + \BV) -\PP(\BY)}_{F}  \le \norm{\BPhi_{0}+\BV- \BY}_{F}.
\]

\begin{lemma} \label{lemma: Lipschitz}
For any $\BV$ in $T_{\BPhi_0}(\St(d,p))$, it holds that
\[
\| R_{\BPhi_0}(\BV) - \BPhi_{0} \|_{F} \leq \|\BV\|_{F}.
\]
Suppose $\|\BPhi-\BPhi_{0}\| < 1$ for some $\BPhi\in\St(d,p)$, then $\BV = R_{\BPhi_{0}}^{-1}(\BPhi)$ exists such that $R_{\BPhi_0}(\BV) = \BPhi$. 
Moreover, if $\| \BPhi - \BPhi_{0} \| \leq 1/8,$ then
\[
\| \BV  \|_{F} \leq 2\|\BPhi - \BPhi_{0}\|.
\]
\end{lemma}
\begin{proof}[\bf Proof of Lemma~\ref{lemma: Lipschitz}]
Lemma~\ref{lem:retraction} implies
\begin{align*}
\| R_{\BPhi_0}(\BV) - \BPhi_{0} \|_F & = \| (\I_d + \BV\BV^{\top})^{-1/2} (\BPhi_{0} + \BV) - \BPhi_{0} \| \leq \|\BPhi_{0} + \BV - \BPhi_{0}\|_F = \|\BV\|_F
\end{align*}
for any $\BV$ on the tangent space.

Next, we will study whether the inverse $R_{\BPhi_0}^{-1}(\BPhi)$ exists for $\BPhi$ sufficiently close to $\BPhi_0$ and whether it is also locally Lipschitz w.r.t. $\BPhi\in\St(d,p)$. 

\noindent {\bf Existence of $R^{-1}_{\BPhi_0}(\cdot)$}. It is easy to see $R^{-1}_{\BPhi_0}(\BPhi_0) = 0$. 
For any $\BV \in T_{\BPhi_{0}}(\St(d,p))$, $\BY : = \BPhi_{0} + \BV$ satisfies
\begin{equation}\label{eq:Rinv}
\BY\BPhi_{0}^{\top} + \BPhi_{0}\BY^{\top} = 2\I,~~~
(\BY \BY^{\top})^{-1/2} \BY = \BPhi.
\end{equation}
We now show that~\eqref{eq:Rinv} has a solution in $\BY$ provided that $\|\BPhi - \BPhi_0\| < 1$.

Let $\BY$ be in the form of $\BY = \BSigma\BPhi$ with $\BSigma\succ 0$, and then it automatically satisfies the second equation in~\eqref{eq:Rinv}. Plugging it into the first constraint in~\eqref{eq:Rinv} gives 
\[
\BSigma\BPhi\BPhi_{0}^{\top} + \BPhi_{0}\BPhi^{\top}\BSigma = 2\I.
\]
The equation has an explicit solution for $\BSigma$ provided that the real part of every eigenvalue of $\BPhi\BPhi_{0}^{\top}$ is positive, i.e., 
\[
\BSigma = 2\int_0^{\infty} e^{-t\BPhi_{0}\BPhi^{\top}} e^{-t \BPhi\BPhi_{0}^{\top}}\diff t
\]
where the integral is finite if $\Re \lambda_i(\BPhi_0\BPhi^{\top}) > 0.$
It is straightforward to verify that $\BY = \BSigma\BPhi$ satisfies~\eqref{eq:Rinv}:
\begin{align*}
\BSigma\BPhi\BPhi_{0}^{\top} + \BPhi_{0}\BPhi^{\top}\BSigma  & = 2 \int_0^{\infty} \left[
e^{-t\BPhi_{0}\BPhi^{\top}} e^{-t \BPhi\BPhi_{0}^{\top}} \BPhi\BPhi_{0}^{\top} + \BPhi_{0}\BPhi^{\top}e^{-t\BPhi_{0}\BPhi^{\top}} e^{-t \BPhi\BPhi_{0}^{\top}}\right] \diff t \\
& = - 2 \int_0^{\infty} \diff [e^{-t \BPhi_{0}\BPhi^{\top}} e^{-t\BPhi\BPhi_{0}^{\top}}] = 2\I_d.
\end{align*}
It suffices to show that $\|\BPhi - \BPhi_0\| < 1$ ensures $\Re \lambda_i(\BPhi_0\BPhi^{\top}) > 0$.
Suppose $\|\BPhi - \BPhi_{0}\| \leq \delta$ for some $\delta < 1$, then
\[
\BPhi\BPhi_{0}^{\top} = \I + (\BPhi-\BPhi_{0})\BPhi_{0}^{\top}.
\]
Note that for any normalized eigenvector $\bu$ of $(\BPhi-\BPhi_{0})\BPhi_{0}^{\top}$ with corresponding eigenvalue $\lambda$, we have
\[
| \bu^* [(\BPhi-\BPhi_{0})\BPhi_{0}^{\top}] \bu| = |\lambda |\leq \|(\BPhi-\BPhi_{0})\BPhi_{0}^{\top}\|.
\]
As a result, if $\|\BPhi - \BPhi_{0}\| \leq \delta$, then all the eigenvalues of $\BPhi\BPhi_{0}^{\top}$ have their real parts at least $1-\delta$, which is of course positive for $\delta<1.$
\vskip0.25cm

\noindent{\bf Lipschitz continuity of $R^{-1}_{\BPhi_0}(\BPhi)$}.
Given the existence of the inverse mapping, we consider upper bounding $R_{\BPhi_0}^{-1}(\BPhi)$ via $\|\BPhi - \BPhi_0\|.$
Let $\BDelta := \BPhi-\BPhi_{0}$ and $\BV = R_{\BPhi_0}^{-1}(\BPhi_{0} + \BDelta)$, i.e.,
\[
(\I_d + \BV\BV^{\top})^{-1/2}(\BPhi_{0}+\BV) = \BPhi_{0}+\BDelta = \BPhi.
\]
Our goal is to control $\|\BV\|$ by $\|\BDelta\|.$
Let
\[
h(t) := (\I_d + t^2\BV\BV^{\top})^{-1/2}(\BPhi_{0}+t\BV)
\]
be a matrix-valued function. It derivative equals
\begin{align*}
h'(t) & = (\I_d + t^2\BV\BV^{\top})^{-1/2}\BV - t(\I_d + t^2\BV\BV^{\top})^{-3/2} \BV\BV^{\top}(\BPhi_{0}+t\BV).
\end{align*}
Then it holds that
\begin{align*}
\BDelta & = \int_0^t h'(s)\diff s  = \int_0^t \left[ (\I_d + s^2\BV\BV^{\top})^{-1/2}\BV - s(\I_d + s^2\BV\BV^{\top})^{-3/2} \BV\BV^{\top}(\BPhi_{0}+s\BV)  \right] \diff s.
\end{align*}

Multiplying both sides by $\BV^{\top}$ gives
\begin{align*}
\BDelta \BV^{\top}&  = \int_0^t \left[ (\I_d + s^2\BV\BV^{\top})^{-1/2}\BV\BV^{\top} - s(\I_d + s^2\BV\BV^{\top})^{-3/2} \BV\BV^{\top}(\BPhi_{0}+s\BV)\BV^{\top}  \right] \diff s
\end{align*}
and taking the Frobenius norm gives
\begin{align*}
\|\BDelta\|\| \BV\|_{F} & \geq \left\| \int_0^1 \left[ (\I_d + s^2\BV\BV^{\top})^{-1/2}\BV\BV^{\top} - s(\I_d + s^2\BV\BV^{\top})^{-3/2} \BV\BV^{\top}(\BPhi_{0}+s\BV)\BV^{\top}  \right] \diff s \right\|_{F} \\
& \geq \left\| \int_0^1 (\I_d + s^2\BV\BV^{\top})^{-1/2}\BV\BV^{\top} \diff s \right\|_{F}- \left\| \int_0^1 s(\I_d + s^2\BV\BV^{\top})^{-3/2} \BV\BV^{\top}(\BPhi_{0}+s\BV)\BV^{\top}  \diff s \right\|_{F}  \\
& \geq (1 + \|\BV\|_{F}^2)^{-1/2} \|\BV\|_{F}^2 - \int_0^1 s \diff s \cdot \|\BV\|_{F}^3 (\|\BV\|_{F} + 1) \\
& = (1 + \|\BV\|_{F}^2)^{-1/2} \|\BV\|_{F}^2 -  \frac{1}{2}  \|\BV\|_{F}^3 (\|\BV\|_{F} + 1)
\end{align*}
where $\|\BV\|_{F} \leq 1$ is assumed. Therefore, we have
\[
(1 + \|\BV\|_{F}^2)^{-1/2} \|\BV\|_{F} -  \frac{1}{2}\|\BV\|_{F}^2(\|\BV\|_{F}+1) \leq \|\BDelta\|.
\]
Define the following function
\[
\ell(x) := \frac{x}{\sqrt{1+x^2}} - \frac{1}{2}x^2(x+1).
\]
We can see that $\ell(\|\BV\|_F) \leq \|\BDelta\|.$
For $0 \leq x \leq 1/4$, it holds
\begin{align*}
\ell'(x) & = \frac{1}{(1+x^2)^{3/2}} - \left(\frac{3x^2}{2} + x\right)  > \frac{1}{2},~~~~~~\ell(0) = 0.
\end{align*}
Therefore, we have $\|\BDelta\| \geq  \|\BV\|_{F}/2$ for $\|\BV\|_F \leq 1/4$ and thus $\|\BV\|_{F} \leq 2\|\BDelta\|$
for $\|\BDelta\| \leq 1/8$.
\end{proof}

\begin{lemma}
For the retraction, $R_{\BPhi_0}(t\BV)$ satisfies
\[
R_{\BPhi_0}(t\BV) = \BPhi_{0} + t\BV  - \frac{t^2}{2}\BV\BV^{\top}\BPhi_{0} - \frac{t^3}{2}\BV\BV^{\top}\BV + \BE
\]
where $\|\BE\| \leq t^4\|\BV\|^4 $ for $\|\BV\|_{F} \leq 1/4$ and $t\leq 1$.
Moreover, we have $g_V(t) = \lag \BA, R_{\BPhi_0}(t\BV)\rag$ satisfies
\begin{equation}\label{eq:gv_taylor}
g_V(t) := \lag \BA, R_{\BPhi_0}(t\BV)\rag  = \lag \BA, \BPhi_{0}\rag + t\lag \BA,\BV\rag - \frac{t^2}{2} \lag \BA\BPhi_{0}^{\top}, \BV\BV^{\top}\rag - \frac{t^3}{2}\lag \BA, \BV\BV^{\top}\BV\rag + E_V
\end{equation}
with residual  $|E_V| \leq t^4\|\BA\|\|\BV\|_{F}^4 $ for $\|\BV\|_{F}\leq 1/4$ and $t\leq 1.$
\end{lemma}

\begin{proof}
Consider $R_{\BPhi_0}(t\BV) = (\I + t^2\BV\BV^{\top})^{-1/2}(\BPhi_{0}+t\BV).$
Define $\BK(t) := (\I + t^2 \BV\BV^{\top})^{-1/2}$. Then it holds that
{
\begin{align*}
\BK'(t) & = -t(\I + t^2 \BV\BV^{\top})^{-3/2}\BV\BV^{\top}, \\
\BK''(t) & = -(\I + t^2 \BV\BV^{\top})^{-3/2}\BV\BV^{\top} + 3t^2 (\I + t^2\BV\BV^{\top})^{-5/2} (\BV\BV^{\top})^2, \\
\BK'''(t) 
& =  9t (\I + t^2\BV\BV^{\top})^{-5/2} (\BV\BV^{\top})^2 - 15t^3 (\I + t^2\BV\BV^{\top})^{-7/2}(\BV\BV^{\top})^3, \\
\BK^{(4)}(t) & = 9 (\I + t^2\BV\BV^{\top})^{-5/2} (\BV\BV^{\top})^2  - 90t^2 (\I + t^2\BV\BV^{\top})^{-7/2}(\BV\BV^{\top})^3 \\
& \qquad + 105t^4(\I + t^2\BV\BV^{\top})^{-9/2}(\BV\BV^{\top})^4.
\end{align*}}
The last term $\BK^{(4)}(t)$ satisfies
\[
-90t^2\|\BV\|^2(\BV\BV^{\top})^2 \preceq \BK^{(4)}(t) \preceq (9 + 105t^4\|\BV\|^4)(\BV\BV^{\top})^2
\]
where $\|(\I + t^2\BV\BV^{\top})^{-1/2}\| \leq 1$ for any $t.$

The Taylor approximation with the integral form of the remainder gives
\begin{align*}
\BK(t)  & = \I - \frac{\BV\BV^{\top} }{2}  t^2 + \frac{1}{6} \int_0^t \BK^{(4)}(s) s^3 \diff s
\end{align*}
and the remainder satisfies
\begin{align*}
& \int_0^t \BK^{(4)}(s) s^3 \diff s \succeq -90\int_0^t s^5 \diff s \cdot \|\BV\|^2 (\BV\BV^{\top})^2 \succeq - 15 \|\BV\|^2 t^6(\BV\BV^{\top})^2, \\
& \int_0^t \BK^{(4)}(s) s^3 \diff s \preceq \left( 9 \int_0^t s^3 \diff s+ 105 \|\BV\|^4 \int_0^t s^7\diff s\right) (\BV\BV^{\top})^2  =\left( \frac{9 t^4}{4}  + \frac{105\|\BV\|^4t^8}{8} \right)(\BV\BV^{\top})^2.
\end{align*}
In particular, if $\|\BV\|_F \leq 1/4$, then
\[
-3t^4(\BV\BV^{\top})^2 \preceq \int_0^t \BK^{(4)}(s) s^3 \diff s  \preceq 3t^4(\BV\BV^{\top})^2,
\]
and thus
\[
-\frac{t^4(\BV\BV^{\top})^2}{2} \preceq \BK(t) - \I + \frac{t^2\BV\BV^{\top}}{2}   \preceq \frac{t^4(\BV\BV^{\top})^2}{2}.
\]

Then
\begin{align*}
R_{\BPhi_0}(t\BV) & : = (\I + t^2\BV\BV^{\top})^{-1/2}(\BPhi_{0}+ t\BV) = \BK(t) (\BPhi_0 + t\BV)\\
& = \left(\I - \frac{\BV\BV^{\top}}{2} t^2 \right) (\BPhi_{0}+ t\BV) +  \left(\BK(t) - \I + \frac{\BV\BV^{\top}}{2} t^2\right) (\BPhi_{0}+ t\BV) \\
& = \BPhi_{0} + t\BV - \frac{\BV\BV^{\top}\BPhi_{0}}{2} t^2 - \frac{\BV\BV^{\top} \BV}{2}t^3 + \underbrace{ \left(\BK(t) - \I +\frac{\BV\BV^{\top}}{2} t^2\right) (\BPhi_{0} + t\BV)}_{\BE}.
\end{align*}
The error term satisfies:
\begin{align*}
\left\| \BE \right\| & \leq 2\left\|\BK(t) - \I +\frac{\BV\BV^{\top}}{2} t^2 \right\|  = \frac{1}{3}\left\| \int_0^t \BK^{(4)}(s) s^3 \diff s \right\|  \leq t^4\|\BV\|^4.
\end{align*}

For $g_V(t)$, it holds
\[
g_{V}(t) = \left\lag \BA, \BPhi_{0} + t\BV - \frac{\BV\BV^{\top}\BPhi_{0}}{2} t^2 - \frac{\BV\BV^{\top} \BV}{2}t^3 \right\rag + \left\lag \BK(t) - \I + \frac{\BV\BV^{\top}}{2} t^2,\BA(\BPhi_{0} + t\BV)^{\top} \right\rag
\]
and the residual satisfies
\begin{align*}
\left| \left\lag \BK(t) - \I + \frac{\BV\BV^{\top}}{2} t^2,\BA(\BPhi_{0} + t\BV)^{\top} \right\rag \right| 
& \leq \| \BA (\BPhi_0 + t\BV)^{\top}\| \left\|  \BK(t) - \I+ \frac{t^2\BV\BV^{\top}}{2}\right\|_* \\
& \leq 2t^4\|\BA\| \cdot \frac{\|\BV\|_F^4}{2} = t^4\|\BA\|\|\BV\|_F^4.
\end{align*}
\end{proof}

The following proposition characterizes the one-step Newton iteration applied $f(\BPhi)$, and its distance to the global minimizer.
\begin{proposition}\label{prop:2nderror}
Consider $f(\BPhi) = \lag \BA, \BPhi\rag$ where $\BPhi\in\St(d,p)$. The followings hold:
\begin{enumerate}
\item The one-step Newton iteration applied to $f(\BPhi)$ at $\BPhi_0$ gives
\[
\widetilde{\BPhi} = R_{\BPhi_0}(\widetilde{\BV}) = \PP(\BPhi_0+\widetilde{\BV}),~~ \widetilde{\BV} = - (\BH^{-1}\BA +\BPhi_{0})
\]
where $\widetilde{\BV}$ minimizes the quadratic approximation of $\lag \BA, R_{\BPhi_0}(\BV)\rag$.
\item Let $\widehat{\BV}$ be the element on the tangent space $T_{\BPhi_0}(\St(d,p))$ such that $\widehat{\BPhi}=R_{\BPhi_0}(\widehat{\BV}) $. Suppose 
\[
\BH = - \frac{1}{2} (\BA\BPhi_{0}^{\top} + \BPhi_{0}\BA^{\top}) \succ 0,~~~ \max\{\|\widehat{\BV}\|_F, \|\widetilde{\BV}\|_F\} \leq \frac{1}{4},
\]
then it holds
\[
\|\widehat{\BPhi} - \widetilde{\BPhi}\|_F \leq \|\widehat{\BV} - \widetilde{\BV}\|_F \leq \frac{3}{2} \left( \frac{\|\BA\|}{\lambda_{\min}(\BH)} +  \sqrt{  \frac{\|\BA\|}{\lambda_{\min}(\BH)}} \right) (\|\widehat{\BV}\|_F^2  + \|\widetilde{\BV}\|_F^2 ).
\]
\end{enumerate}

\end{proposition}
\begin{proof}[\bf Proof of Proposition~\ref{prop:2nderror}]
The existence of $\widehat{\BV}$ on $T_{\BPhi_0}(\St(d,p))$ such that $\widehat{\BPhi} = R_{\BPhi_0}(\widehat{\BV})$ is guaranteed by Lemma~\ref{lemma: Lipschitz}, provided that $\|\widehat{\BPhi} - \BPhi_0\| < 1/8.$ Moreover, $\|\widehat{\BV}\|_F \leq 2\|\BPhi - \BPhi_0\| \leq 1/4.$
For $g_V(1)$ that is defined in~\eqref{eq:gv_taylor}, we consider its second-order approximation
\begin{align*}
\overline{g}_V(1) & = \lag \BA,\BPhi_{0}\rag  + \lag \BA,\BV\rag + \frac{1}{2} \lag \BH, \BV\BV^{\top} \rag \\
&  = \lag \BA,\BPhi_{0}\rag + \left\lag \BA -  \frac{\BA\BPhi_{0}^{\top} +\BPhi_{0}\BA^{\top}}{2}\BPhi_{0},\BV\right\rag + \frac{1}{2} \lag \BH, \BV\BV^{\top} \rag 
\end{align*}
where $\BV$ is on the tangent space, i.e., $\BPhi_{0}\BV^{\top} + \BV\BPhi_{0}^{\top} = 0$ and $\BH = -(\BA\BPhi_0^{\top} + \BPhi_0\BA^{\top})/2.$

It is straightforward to verify that 
\begin{align*}
\widetilde{\BV} & := - \BH^{-1}\left( \BA -  \frac{\BA\BPhi_{0}^{\top} + \BPhi_{0}\BA^{\top}}{2}\BPhi_{0}\right) = - \BH^{-1}\left( \BA + \BH\BPhi_{0}\right) = -(\BH^{-1}\BA + \BPhi_{0})
\end{align*}
is in $T_{\BPhi_0}(\St(d,p))$ that minimizes $\bar{g}_V(1)$ over $V$. More precisely, we have
\begin{equation} \label{eqn: newton on tangent}
	\begin{aligned}
\widetilde{\BV}\BPhi_{0}^{\top} + \BPhi_{0}\widetilde{\BV}^{\top} & = - (\BH^{-1}\BA\BPhi_{0}^{\top} + \BPhi_{0}\BA^{\top}\BH^{-1}) - 2\I_d \\
& = - \BH^{-1} (\BA\BPhi_{0}^{\top}\BH + \BH\BPhi_{0}\BA^{\top})\BH^{-1} - 2\I_d \\
& = \frac{1}{2}\BH^{-1} (\BA\BPhi_{0}^{\top} (\BA\BPhi_{0}^{\top}+\BPhi_{0}\BA^{\top}) + (\BA\BPhi_{0}^{\top}+\BPhi_{0}\BA^{\top})  \BPhi_{0}\BA^{\top} ) \BH^{-1}-2\I_d \\
& = \frac{1}{2}\cdot 4\BH^{-1} \BH^2 \BH^{-1}-2\I_d = 0.
\end{aligned}
\end{equation}
Note that Proposition~\ref{prop:2nderror} implies
\begin{align*}
g_{V}(1) & = \lag \BA, \BPhi_{0}\rag + \lag \BA + \BH \BPhi_{0},\BV\rag + \frac{1}{2} \lag \BH, \BV\BV^{\top}\rag - \frac{1}{2}\lag \BA, \BV\BV^{\top}\BV\rag+ E_V 
\end{align*}
where $|E_V| \leq \|\BA\| \|\BV\|^4_F.$

From the equation above, the global optimality of $\widehat{\BV}$ in $g_V(1)$ implies
\begin{equation}\label{eq:gvv}
\overline{g}_{\widehat{V}}(1) - \frac{1}{2}\lag \BA, \widehat{\BV}\widehat{\BV}^{\top}\widehat{\BV}\rag + E_{\widehat{V}} = g_{\widehat{V}}(1) \leq g_{\widetilde{V}}(1) = \overline{g}_{\widetilde{V}}(1) - \frac{1}{2}\lag \BA, \widetilde{\BV}\widetilde{\BV}^{\top}\widetilde{\BV}\rag + E_{\widetilde{V}}
\end{equation}
where $\overline{g}_V(1) = \lag \BA, \BPhi_{0}\rag + \lag \BA + \BH \BPhi_{0},\BV\rag + \lag \BH, \BV\BV^{\top}\rag/2.$
Now it holds that
\begin{align*}
\overline{g}_{\widehat{V}}(1) - \overline{g}_{\widetilde{V}}(1) & = \frac{1}{2} \lag \BH, \widehat{\BV}\widehat{\BV}^{\top}\rag - \frac{1}{2} \lag \BH, \widetilde{\BV}\widetilde{\BV}^{\top}\rag + \lag \BA+\BH\BPhi_{0}, \widehat{\BV} - \widetilde{\BV}\rag \\
& = \frac{1}{2} \lag \BH, \widehat{\BV}\widehat{\BV}^{\top}\rag - \frac{1}{2} \lag \BH, \widetilde{\BV}\widetilde{\BV}^{\top}\rag - \lag \BH\widetilde{\BV}, \widehat{\BV} - \widetilde{\BV}\rag \\
& = \frac{1}{2} \lag \BH, (\widehat{\BV} - \widetilde{\BV})(\widehat{\BV}-\widetilde{\BV})^{\top}\rag \\
& \leq \frac{1}{2}\lag \BA,  \widehat{\BV}\widehat{\BV}^{\top}\widehat{\BV} -\widetilde{\BV}\widetilde{\BV}^{\top}\widetilde{\BV}\rag + (E_{\widetilde{V}} - E_{\widehat{V}})
\end{align*}
where $\BH\widetilde{\BV} = -\BA -  \BH\BPhi_{0}.$
Therefore, we have from~\eqref{eq:gvv} that
\begin{align*}
\lambda_{\min}(\BH)\|\widehat{\BV}-\widetilde{\BV}\|_F^2 & \leq \lag \BH, (\widehat{\BV} - \widetilde{\BV})(\widehat{\BV}-\widetilde{\BV})^{\top}\rag \\
 & \leq \lag \BA, \widehat{\BV}\widehat{\BV}^{\top}\widehat{\BV} -\widetilde{\BV}\widetilde{\BV}^{\top}\widetilde{\BV} \rag+ 2(E_{\widetilde{V}} - E_{\widehat{V}}) \\
 & \leq \|\BA\|\left( \|\widehat{\BV}\|_F^2 + \|\widehat{\BV}\|_F\|\widetilde{\BV}\|_F + \|\widetilde{\BV}\|_F^2\right) \|\widehat{\BV} - \widetilde{\BV}\|_F + 2\| \BA \| (\|\widehat{\BV}\|_F^4 + \|\widetilde{\BV}\|_F^4)
\end{align*}
where
\begin{align*}
\lag \BA, \widehat{\BV}\widehat{\BV}^{\top}\widehat{\BV} -\widetilde{\BV}\widetilde{\BV}^{\top}\widetilde{\BV} \rag 
& = \lag \BA, \widehat{\BV}\widehat{\BV}^{\top}(\widehat{\BV} -\widetilde{\BV}) 
+ (\widehat{\BV}- \widetilde{\BV})\widehat{\BV}^{\top}\widetilde{\BV} 
+ \widetilde{\BV}(\widehat{\BV}-\widetilde{\BV})^{\top}\widetilde{\BV} \rag \\
& \leq \|\BA\|\left( \|\widehat{\BV}\|_F^2 + \|\widehat{\BV}\|_F\|\widetilde{\BV}\|_F + \|\widetilde{\BV}\|_F^2\right) \|\widehat{\BV} - \widetilde{\BV}\|_F. 
\end{align*}
Solving the inequality above leads to
\begin{align*}
\|\widehat{\BV} - \widetilde{\BV}\|_F 
& \leq \frac{\|\BA\|}{\lambda_{\min}(\BH)}\left( \|\widehat{\BV}\|_F^2 + \|\widehat{\BV}\|_F\|\widetilde{\BV}\|_F + \|\widetilde{\BV}\|_F^2\right) + \sqrt{ \frac{2\|\BA\|}{\lambda_{\min}(\BH)} (\|\widehat{\BV}\|_F^4 + \|\widetilde{\BV}\|_F^4)} \\
& \leq \frac{3}{2}\left( \frac{\|\BA\|}{\lambda_{\min}(\BH)} +  \sqrt{  \frac{\|\BA\|}{\lambda_{\min}(\BH)}} \right) (\|\widehat{\BV}\|_F^2  + \|\widetilde{\BV}\|_F^2 )
\end{align*}
where the estimation of $E_{\widetilde{V}}$ and $E_{\widehat{V}}$ follows from~\eqref{eq:gv_taylor}.
\end{proof}

\section{Proof of Theorem~\ref{thm:main}}
\label{ss:propkey}
Before moving to the proof of Theorem~\ref{thm:main}, we need to collect a few useful lemmas from~\cite{L22}. Under the conditions of Theorem~\ref{thm:main}, the MLE $\widehat{\BZ}\in \Od(d)^{\otimes}$ satisfies with probability at least $1-n^{-2}$ that
\begin{equation}\label{eq:assumption}
\begin{aligned}
& \min_{\BQ\in\Od(d)} \max_{1\leq i\leq n} \|\widehat{\BZ}_i - \BZ_i\BQ\|_F \lesssim  \SNR^{-1}(\sqrt{d} + \sqrt{\log n}),\\
& \min_{\BQ\in\Od(d)}\|\widehat{\BZ} - \BZ\BQ\|_F \lesssim   \SNR^{-1}(\sqrt{d} + \sqrt{\log n}) \cdot \sqrt{n},
\end{aligned}
\end{equation}
where $\BW_i$ is the $i$-th block column of $\BW$. The Gaussian random matrix $\BW$ satisfies $\norm{\BW} \le 3\sqrt{nd}$
with probability at least $1-\exp(-nd)$. 

Throughout our discussion, we assume $\BQ= \PP(\BZ^{\top}\widehat{\BZ})=\I_d$ without loss of generality as $\widehat{\BZ}\BQ^{\top}(\widehat{\BZ}\BQ^{\top})^{\top} = \widehat{\BZ}\widehat{\BZ}^{\top}$ for any $\BQ \in \Od(d)$. As discussed in Section~\ref{ss:proof sketch mle}, we will consider
\begin{equation}\label{def:Amle}
f(\BPhi) = \lag \BA, \BPhi\rag,~~\BPhi\in\Od(d),~~~~\BA = -\sum_{j=1}^n \BC_{ij}\widehat{\BZ}_j = - (\BZ_i \BZ^{\top}\widehat{\BZ} + \sigma\BW_i^{\top}\widehat{\BZ}).
\end{equation}

\subsection{Supporting Lemmas}

\begin{lemma}\label{lem:zs}
For  any $\BS$ satisfying $\BS^{\top}\BS = n\BI_{d}$, it holds that
\[
\|\BZ^{\top}\BS - n\BQ\|_F \le  \frac{1}{2}\|\BS - \BZ\BQ\|_F^2
\]
where $\BQ = \PP(\BZ^{\top}\BS)$ is the polar retraction of $\BZ^{\top}\BS.$
\end{lemma}
\begin{proof}[\bf Proof of Lemma~\ref{lem:zs}]
Note that
\begin{align*}
\| \BZ^{\top}\BS - n\BQ \|_F^2 & = \|\BZ^{\top}\BS\|_F^2 -2n\|\BZ^{\top}\BS\|_* + n^2 d \\
&  = \sum_{i=1}^d (n - \sigma_i(\BZ^{\top}\BS))^2 \leq (nd - \|\BZ^{\top}\BS\|_*)^2 \leq \frac{1}{4} \|\BS - \BZ\BQ\|_F^4
\end{align*}
where
$\|\BS - \BZ\BQ\|_F^2 = 2nd - 2\|\BZ^{\top}\BS\|_*$ and $\sigma_i(\BZ^{\top}\BS)\leq n.$
\end{proof}

\begin{lemma}\label{lem:ws}
Under the conditions of Theorem~\ref{thm:main}, $\BW_i^{\top}\widehat{\BZ}$ satisfies
\[
\BW_i^{\top}\widehat{\BZ} \overset{d}{=} \sqrt{n-1} \BG_i + O(\sqrt{nd} \SNR^{-1}(\sqrt{d} + \sqrt{\log n})),~~1\leq i\leq n,
\]
where $\BG_i$ is a $d\times d$ random matrix with i.i.d. standard Gaussian entries.
\end{lemma}

\begin{proof}[\bf Proof of Lemma~\ref{lem:ws}]
Note that $\BW_i$ and $\widehat{\BZ}$ are statistically dependent. The idea is to approximate $\widehat{\BZ}$ by an auxiliary approximation via the leave-one-out technique.
Consider
\begin{equation}\label{def:Ci}
\BC^{(i)} = \BZ\BZ^{\top} + \sigma \BW^{(i)}, \qquad
\BW^{(i)}_{jk} = 
\begin{cases}
\BW_{jk}, & \text{ if }j \neq i \text{ and }k\neq i, \\
0, & \text{ if } j = i \text{ or } k = i.
\end{cases}
\end{equation}
In other words, $\BW^{(i)}$ equals $\BW$ except the $i$-th block row and column. For the generalized power method (GPM) applied to $\BC^{(i)}$:
\[
\widehat{\BZ}^{\ell,i} =\PP_n(\BC^{(i)} \widehat{\BZ}^{\ell-1,i}),~~1\leq i\leq n,~~~1\leq \ell\leq L
\]
where $\PP_n(\cdot)$ is the retraction on $\Od(d)^{\otimes n}.$ In particular, $\widehat{\BZ}^{\ell}$ is the iterate generated by the GPM applied to $\BC$. The GPM is defined as: for $t=0,1,2,\ldots$, $\BZ^{t+1} = \mathcal{P}_{n}(\BC\BZ^{t})$ where $\BZ^{0}$ is the top-$d$ eigenvectors of $\BC$.
As a result, $\widehat{\BZ}^{\ell,i}$ is independent of $\BW_i$ since the iteration does not involve $\BW_i$ where $\BW_i$ is the $i$-th block column of $\BW$ for $1\leq\ell\leq L$.

From~\cite[Lemma 4.13 and 4.14]{L22}, it holds
\begin{align*}
& d_F(\widehat{\BZ}^{\ell},\widehat{\BZ}^{\ell,i}) \leq \eps\sqrt{d},~~~1\leq i\leq n,~~1\leq \ell\leq L, \\
& d_F(\widehat{\BZ}, \widehat{\BZ}^L) \leq 2^{-L+1}\eps\sqrt{nd},~~L = n
\end{align*}
with probability at least $1-n^{-2}$, provided that
$\sigma \leq c_0\eps\sqrt{n}/(\sqrt{d} + \sqrt{\log n})$
for some constant $c_0>0.$
In particular, we set
\[
\sigma = \frac{c_0\eps\sqrt{n}}{\sqrt{d} + \sqrt{\log n}} ~~\Longleftrightarrow~~ \eps\sqrt{d} = c_0^{-1}\SNR^{-1}(\sqrt{d} + \sqrt{\log n}).
\]
As a result, we have
\[
d_F(\widehat{\BZ}, \widehat{\BZ}^{L,i}) \leq d_F(\widehat{\BZ},\widehat{\BZ}^L)+ d_F(\widehat{\BZ}^{\ell},\widehat{\BZ}^{\ell,i})  \leq 2\eps\sqrt{d} \leq 2c_0^{-1}\SNR^{-1}(\sqrt{d} + \sqrt{\log n}),~~\forall 1\leq i\leq n
\]
where $2^{-n+1}\sqrt{n} \leq 1$ for $n\geq 3$ and $L = n.$

Now we consider
\begin{align*}
\BW_i^{\top}\widehat{\BZ} & = \BW_i^{\top}(\widehat{\BZ} - \widehat{\BZ}^{L,i}\BQ^i) + \BW_i^{\top} \widehat{\BZ}^{L,i}(\BQ^i - \I_d) + \BW_i^{\top} \widehat{\BZ}^{L,i}  = : T_1' + T_2' + T_3'
\end{align*}
where $\BQ^i = \PP((\widehat{\BZ}^{L,i})^{\top}\widehat{\BZ} )$ is statistically dependent on $\widehat{\BZ}$ and $\BW_i.$

For the first term $T_1'$, we note that
\[
\|T_1'\|_F \leq \|\BW_i\|\cdot d_F(\widehat{\BZ}, \widehat{\BZ}^{L,i}) \lesssim \sqrt{nd} \SNR^{-1}(\sqrt{d} + \sqrt{\log n})
\]
where $\|\BW_i\|\leq 3\sqrt{nd}$ holds with high probability.
For $T_2'$, it holds that
\begin{align*}
\|T_2'\|_F & \leq \|\BW_i\|\|\widehat{\BZ}\|\|\I_d - \BQ^{i}\|_F \leq \sqrt{nd} \cdot \sqrt{n}\cdot n^{-1/2} \SNR^{-1}(\sqrt{d} + \sqrt{\log n}) \\
& \lesssim \sqrt{nd} \SNR^{-1}(\sqrt{d} + \sqrt{\log n})
\end{align*}
where we use the claim that 
$\|\I_d - \BQ^{i}\|_F \lesssim n^{-1/2} \SNR^{-1}(\sqrt{d} + \sqrt{\log n})$
which will be proven later.
For the third term, we have
$T_3' \overset{d}{=} \sqrt{n-1} \BG_i$
where $\BG_i$ is a $d\times d$ Gaussian random matrix since $\widehat{\BZ}^{L,i}$ is in $\Od(d)^{\otimes n}$ and independent of $\BW_i.$ Therefore, we have the result.
Finally, it remains to prove the claim regarding $\|\BQ^i - \I_d\|_F$.
Without loss of generality, we set
\[
\I_d = \PP(\BZ^{\top}\widehat{\BZ}),~~~\I_d = \PP(\BZ^{\top}\widehat{\BZ}^{L,i}),~~~\BQ^i = \PP((\widehat{\BZ}^{L,i})^{\top}\widehat{\BZ} )
\]
and then it holds
\[
 \|\BZ^{\top}(\widehat{\BZ} - \widehat{\BZ}^{L,i}\BQ^i) \|_F \leq \sqrt{n}\| \widehat{\BZ} - \widehat{\BZ}^{L,i}\BQ^i \|_F \leq 2\eps\sqrt{nd}.
\]

Note that
\begin{equation}\label{eq:QminusId}
\begin{aligned}
\|\I_d - \BQ^i\|_F 
& = \|\PP(\BZ^{\top}\widehat{\BZ}) - \PP(\BZ^{\top}\widehat{\BZ}^{L,i}) \BQ^i \|_F  \leq \frac{2\| \BZ^{\top}(\widehat{\BZ} - \widehat{\BZ}^{L,i}\BQ^i)  \|_F}{\sigma_{\min}(\BZ^{\top}\widehat{\BZ}) + \sigma_{\min}(\BZ^{\top}\widehat{\BZ}^{L,i}) }  \\
& \lesssim \frac{\| \widehat{\BZ} - \widehat{\BZ}^{L,i}\BQ^i  \| }{\sqrt{n}}\lesssim  \eps\sqrt{\frac{d}{n}}  \lesssim n^{-1/2} \SNR^{-1}(\sqrt{d} + \sqrt{\log n})
\end{aligned}
\end{equation}
where the first inequality follows from~\cite{L95} and~\cite[Lemma 4.6]{L22} and Lemma~\ref{lem:zs} gives
\[
\max \{n - \sigma_{\min}(\BZ^{\top}\widehat{\BZ}), n-\sigma_{\min}(\BZ^{\top}\widehat{\BZ}^{L,i})\} \lesssim n (\SNR^{-1}(\sqrt{d} + \sqrt{\log n}))^2.
\]
\end{proof}

\begin{lemma}\label{lem:HessMLE}
Under the conditions of Theorem~\ref{thm:main}, $\BH_i$ defined in~\eqref{def:HI} satisfies with probability at least $1-n^{-2}$ that
	\[
	\|\BH_i - n\I_d\| \lesssim n\SNR^{-1}(\sqrt{d} + \sqrt{\log n}),~~\forall 1\leq i\leq n.
	\]
\end{lemma}

\begin{proof}[\bf Proof of Lemma~\ref{lem:HessMLE}]

By definition, we have
\[
\BH_i = \frac{1}{2} \sum_{j=1}^n \left[ \BZ_i (\BZ^{\top}\widehat{\BZ} + \widehat{\BZ}^{\top}\BZ)\BZ_i^{\top} + \frac{\sigma}{2} (\BW_i^{\top}\widehat{\BZ} \BZ_i^{\top} + \BZ_i \widehat{\BZ}^{\top}\BW_i)\right]
\]
where $\BW_i$ is the $i$th block column of $\BW$ and $\I_d = \PP(\BZ^{\top}\widehat{\BZ}).$
Using Lemma~\ref{lem:zs} and~\eqref{eq:assumption} gives $\| \BZ^{\top}\widehat{\BZ} - n\I_d \|_F \leq \frac{1}{2}\|\widehat{\BZ} - \BZ\|_F^2 \lesssim n (\SNR^{-1}(\sqrt{d} +\sqrt{\log n}))^2.$ Then
\begin{align*}
\|\BH_i - n\I_d\| & \leq  \frac{1}{2} \left\|  \BZ^{\top} \widehat{\BZ} + \widehat{\BZ}^{\top}\BZ - 2n\I_d\right\|  + \frac{\sigma}{2} \left\| \BW_i^{\top}\widehat{\BZ} \BZ_i^{\top} + \BZ_i \widehat{\BZ}^{\top}\BW_i\right\| \\
& \leq \| \BZ^{\top}\widehat{\BZ} - n\I_d \|_F + \sigma \|\BW_i^{\top}\widehat{\BZ}\|_F \\
& \lesssim n (\SNR^{-1}(\sqrt{d} +\sqrt{\log n}))^2 + n\SNR^{-1}(\sqrt{d} +\sqrt{\log n})
\end{align*}
where the estimation of $\|\BW_i^{\top}\widehat{\BZ}\|_F$ follows from Lemma~\ref{lem:ws}.
Therefore, $\|\BH_i - n\I_d\|\lesssim n\SNR^{-1}(\sqrt{d} + \sqrt{\log n}).$
\end{proof}


\subsection{Proof of Proposition~\ref{prop:tildeV}}
\begin{proof}[\bf Proof of Proposition~\ref{prop:tildeV}]

Note that $\|\widehat{\BV}_i - \widetilde{\BV}_i\|_F\lesssim (\SNR^{-1}(\sqrt{d} + \sqrt{\log n}))^2$, and then it suffices to find the distribution of $\widehat{\BV}_i$, which can be decomposed into three parts:
\begin{align*}
- \widetilde{\BV}_i  
& = \frac{1}{2}\BH_i^{-1} ( \BA\BZ_i^{\top} - \BZ_i\BA^{\top} )\BZ_i \\
& = \frac{1}{2}(\BH_i^{-1} - n^{-1}\I_d) ( \BA\BZ_i^{\top} - \BZ_i\BA^{\top} )\BZ_i  + \frac{1}{2n} ( \BA\BZ_i^{\top} - \BZ_i\BA^{\top} )\BZ_i \\
& = \frac{1}{2}(\BH_i^{-1} - n^{-1}\I_d) ( \BA\BZ_i^{\top} - \BZ_i\BA^{\top} )\BZ_i  + \frac{1}{2n} \BZ_i (\BZ^{\top}\widehat{\BZ} - \widehat{\BZ}^{\top}\BZ) \\
& \qquad + \frac{\sigma}{2n} (\BW_i^{\top}\widehat{\BZ}\BZ_i^{\top} - \BZ_i \widehat{\BZ}^{\top}\BW_i) \BZ_i  = : T_1 + T_2 + T_3
\end{align*}
where $\BA$ is defined in~\eqref{def:Amle} and
\begin{align*}
\BA\BZ_i^{\top} - \BZ_i\BA^{\top} & = (\BZ_i\BZ^{\top}\widehat{\BZ} + \sigma\BW_i^{\top}\widehat{\BZ}) \BZ_i^{\top} - \BZ_i (\BZ_i\BZ^{\top}\widehat{\BZ} + \sigma\BW_i^{\top}\widehat{\BZ})^{\top} \\
& = \BZ_i (\BZ^{\top}\widehat{\BZ} - \widehat{\BZ}^{\top}\BZ)\BZ_i^{\top} + \sigma (\BW_i^{\top}\widehat{\BZ}\BZ_i^{\top} - \BZ_i \widehat{\BZ}^{\top}\BW_i).
\end{align*}

We now estimate each term above: for the first term, note that
\begin{align*}
\| \BH_i^{-1} - n^{-1}\I_d\| & = n^{-1} \|\BH_i^{-1}\| \|\BH_i - n\I_d\| \lesssim n^{-1}\SNR^{-1}(\sqrt{d} + \sqrt{\log n})
\end{align*}
and Lemma~\ref{lem:zs} gives
\[
\|\BZ^{\top}\widehat{\BZ} - n\I_d \|_F \le  \frac{1}{2}\|\widehat{\BZ} - \BZ\|_F^2 \lesssim n (\SNR^{-1}(\sqrt{d} + \sqrt{\log n}))^2.
\]
As a result, it holds $\|\BZ^{\top}\widehat{\BZ} - \widehat{\BZ}^{\top}\BZ\|_F\lesssim  n (\SNR^{-1}(\sqrt{d} + \sqrt{\log n}))^2$ and
\begin{align*}
\|\BA\BZ_i^{\top} - \BZ_i\BA^{\top}\|_F & \leq \|\BZ^{\top}\widehat{\BZ} - \widehat{\BZ}^{\top}\BZ\|_F + 2\sigma \|\BW_i^{\top}\widehat{\BZ}\|_F  \lesssim n\SNR^{-1} (\sqrt{d} + \sqrt{\log n}). 
\end{align*}

Therefore, the first term $T_1$ is upper bounded via
\begin{align*}
\|T_1\|_F & \leq \frac{1}{2}\|(\BH_i^{-1} - n^{-1}\I_d) ( \BA\BZ_i^{\top} - \BZ_i\BA^{\top} )\BZ_i \|_F \\
& \leq \|\BH_i^{-1} - n^{-1}\I_d\| \|\BA\BZ_i^{\top} - \BZ_i\BA^{\top}\|_F  \lesssim (\SNR^{-1} (\sqrt{d} + \sqrt{\log n}))^2
\end{align*}
and the second term is bounded via
\[
\|T_2\|_F \leq \frac{1}{2n} \| \BZ_i (\BZ^{\top}\widehat{\BZ} - \widehat{\BZ}^{\top}\BZ)\|_F \lesssim  (\SNR^{-1}(\sqrt{d} + \sqrt{\log n}))^2.
\]
It remains to estimate the third term $T_3$: from Lemma~\ref{lem:ws}, we note that
\[
\BW_i^{\top}\widehat{\BZ}\BZ_i^{\top} - \BZ_i\widehat{\BZ}^{\top}\BW_i = \sqrt{n-1}(\BG_i - \BG_i^{\top}) + O(\sqrt{nd} \SNR^{-1}(\sqrt{d} + \sqrt{\log n}))
\]
where $\BG_i$ is a Gaussian random matrix and so is $\BG_i\BZ_i^{\top}$. As a result, we have
\begin{align*}
T_3 &  \overset{d}{=} \frac{\sigma\sqrt{n-1}}{2n}(\BG_i - \BG_i^{\top})\BZ_i + O(\SNR^{-1}(\sqrt{d} + \sqrt{\log n}))^2.
\end{align*}
Combining the estimation of $T_1$, $T_2$, and $T_3$ yields
\[
\widetilde{\BV}_i \overset{d}{=} \frac{\sigma\sqrt{n-1}}{2n}(\BG_i - \BG_i^{\top})\BZ_i  + O(\SNR^{-1}(\sqrt{d} + \sqrt{\log n}))^2
\]
since $T_3$ is the leading term. As a by-product, we have
$\|\widetilde{\BV}_i\|_F \lesssim \SNR^{-1}(\sqrt{d} +\sqrt{\log n})$ for all $1\leq i\leq n$.
\end{proof}

\subsection{Proof of Proposition~\ref{prop:newtonstep}}

\begin{proof}[\bf Proof of Proposition~\ref{prop:newtonstep}]
Let $\widehat{\BV}_i$ and $\widetilde{\BV}_{i}$ be two tangent vectors on $\Od(d)$ at $\BZ_i$ such that 
\[
R_{\BZ_i}(\widehat{\BV}_i) = \widehat{\BZ}_i,~~~~ R_{\BZ_i}(\widetilde{\BV}_i) = \widetilde{\BZ}_i
\]
where $R_{\BZ_i}(\cdot) = \PP(\BZ_i + \cdot)$ is defined in~\eqref{eqn: polar retraction}.
The existence of $\widehat{\BV}_i$ is guaranteed if $\|\widehat{\BZ}_i - \BZ_i\| < 1$, which is verified in Lemma~\ref{lemma: Lipschitz}.

The idea is to apply Lemma~\ref{prop:2nderror} with $\BPhi_0 = \BZ_i$ and $\BA$ in~\eqref{def:Amle}:
\begin{align} \label{eqn: lemma311main}
	\|\widehat{\BZ}_{i} - \widetilde{\BZ}_{i}\|_F 
	& \leq \frac{3}{2} \left( \frac{\|\sum_{j=1}^n \BC_{ij}\widehat{\BZ}_j\|}{\lambda_{\min}(\BH_i)} +  \sqrt{  \frac{\|\sum_{j=1}^n \BC_{ij}\widehat{\BZ}_j\|}{\lambda_{\min}(\BH_i)}} \right) (\|\widehat{\BV}_i\|_F^2  + \|\widetilde{\BV}_i\|_F^2 ).
\end{align}
Therefore, we need to upper bound $\norm{\sum_{j=1}^n\BC_{ij}\widehat{\BZ}_{j}}$, $\|\widehat{\BV}_i\|_F$, and $\|\widetilde{\BV}_i\|_F$:
\begin{equation}\label{eqn: spectral norm of A}
\begin{aligned}
	\norm{\sum_{j=1}^n\BC_{ij}\widehat{\BZ}_{j}} 
	&\le \norm{ \BZ_{i}\BZ^{\top}\widehat{\BZ}} + \sigma\norm{\BW_i^{\top}\widehat{\BZ}}_{F} \lesssim n + \SNR^{-1}(\sqrt{d}+\sqrt{\log n}) \cdot n \lesssim n.
\end{aligned}	
\end{equation}
Using the Lipschitz continuity of $R_{\BZ_{i}}^{-1}$ (Lemma~\ref{lemma: Lipschitz}), $\|\widehat{\BV}_i\|_F$ is upper bounded via 
\begin{align} \label{eqn: F_norm of V_hat}
	\norm{\widehat{\BV}_{i}}_{F} \le 2 \norm{\widehat{\BZ}_{i} - \BZ_{i}}_{F} \lesssim \SNR^{-1}(\sqrt{d}+\sqrt{\log n}).
\end{align}
The bound on $\|\widetilde{\BV}_i\|_F$ satisfies the same bound $\|\widetilde{\BV}_i\|_F\lesssim \SNR^{-1}(\sqrt{d}+\sqrt{\log n})$, following from Proposition~\ref{prop:tildeV}.

According to Lemma~\ref{lem:HessMLE}, we have $\|\BH_i - n\I_d\|\lesssim n\SNR^{-1}(\sqrt{d} + \sqrt{\log n})$ which implies $\lambda_{\min}(\BH_i) \geq n(1 - c\SNR^{-1}(\sqrt{d} +\sqrt{\log n}))$ for some $c>0.$ 
Thus
\[
 \frac{\|\sum_{j=1}^n \BC_{ij}\widehat{\BZ}_j\|}{\lambda_{\min}(\BH_i)} \leq C
\]
for some constant scalar $C>0$. Plugging the estimations of $\widetilde{\BV}_i$ and $\widehat{\BV}_i$ into~\eqref{eqn: lemma311main}, we have:
\[
\|\widehat{\BZ}_{i} - \widetilde{\BZ}_{i}\|_F \leq \|\widehat{\BV}_{i} - \widetilde{\BV}_{i}\|_F \lesssim \|\widehat{\BV}_i\|_F^2  + \|\widetilde{\BV}_i\|_F^2 \lesssim (\SNR^{-1}(\sqrt{d} + \sqrt{\log n}))^2.
\]
\end{proof}

\section{Proof of Theorem~\ref{thm:spectral_main}}
Define
\begin{equation}\label{def:Aisp}
\BA_{i} := -  \sum_{j=1}^n \BC_{ij}\BSeig_j = - ( \BZ_i\BZ^{\top}\BSeig + \sigma\BW_i^{\top}\BSeig)
\end{equation}
where $\widehat{\BS}^{\eig}$ is the global minimizer to~\eqref{eq:spectral} and $\widehat{\BS}_i^{\eig}$ is the $i$-th $d\times d$ block of $\widehat{\BS}^{\eig}.$
First we have the following claim proved.
\begin{lemma}\label{eq:spr_local}
For the $\BA_i$ defined in~\eqref{def:Aisp}, we have
\begin{align*}
\widehat{\BS}_i^{\eig} & = \argmin_{\BS_i^{\top}\BS_i = \BPi_i} \lag \BA_i, \BS_i\rag, ~~~~~~\widehat{\BZ}_i^{\eig}  : = \PP(\widehat{\BS}_i^{\eig}) = \argmin_{\BPhi \in\Od(d)}\lag \BA_i\BPi_i^{1/2}, \BPhi \rag.
\end{align*}
\end{lemma}
\begin{proof}[\bf Proof of Lemma~\ref{eq:spr_local}]
Suppose all $\widehat{\BS}_j^{\eig}$ are fixed for $j\neq i$ in~\eqref{eq:spectral}, then the resulting objective function equals
\[
-2\sum_{j\neq i} \lag \BC_{ij}\widehat{\BS}_j^{\eig}, \BS_i \rag - \lag \BC_{ii}, \BS_i\BS_i^{\top}\rag = -2\sum_{j\neq i} \lag \BC_{ij}\widehat{\BS}_j^{\eig}, \BS_i \rag -\Tr (\BPi_i)
\]
where $\BC_{ii}=\I_d$, $\BS_i^{\top}\BS_i = \BPi_i$ and $\BPi_i = \widehat{\BS}_i^{\eig\top}\widehat{\BS}_i^{\eig}$. The global optimality of $\BS^{\eig}$ implies 
\[
\widehat{\BS}_i^{\eig} = \argmin_{\BS_i^{\top}\BS_i = \BPi_i}-\sum_{j\neq i}\lag \BC_{ij}\widehat{\BS}_j^{\eig}, \BS_i\rag.
\]
By letting $\BPhi  = \BS_i\BPi_i^{-1/2}$, we know that $\BPhi \in\Od(d)$ and then
\begin{align*}
 \min_{\BS_i^{\top}\BS_i = \BPi_i}-\lag \BS_i, \widehat{\BS}_i^{\eig}\rag = \min_{\BS_i^{\top}\BS_i = \BPi_i}-\lag \BS_i \BPi_i^{-1/2}, \widehat{\BS}_i^{\eig} \BPi_i^{1/2}\rag = \min_{\BPhi \in\Od(d)} - \lag \BPhi ,\widehat{\BS}_i^{\eig} \BPi_i^{1/2} \rag.
\end{align*}
Therefore, we have
\begin{align*}
  &\argmin_{\BPhi \in\Od(d)} - \lag \BPhi ,\widehat{\BS}_i^{\eig} \BPi_i^{1/2} \rag = \PP(\widehat{\BS}_i^{\eig} \BPi_i^{1/2}) = (\widehat{\BS}_i^{\eig} \BPi_i\widehat{\BS}_i^{\eig\top})^{-1/2} \widehat{\BS}_i^{\eig} \BPi_i^{1/2} \\
 & ~~~~= (\widehat{\BS}_i^{\eig}\widehat{\BS}_i^{\eig\top})^{-1}\widehat{\BS}_i^{\eig} \BPi_i^{1/2} = \widehat{\BS}_i^{\eig} (\widehat{\BS}_i^{\eig\top}\widehat{\BS}_i^{\eig})^{-1}\BPi_i^{1/2} = \widehat{\BS}_i^{\eig} (\widehat{\BS}_i^{\eig\top}\widehat{\BS}_i^{\eig})^{-1/2} = \PP(\widehat{\BS}_i^{\eig})
\end{align*}
where $\PP(\cdot)$ follows from~\eqref{def:polar} and $\BPi_i = \widehat{\BS}_i^{\eig\top}\widehat{\BS}_i^{\eig}$. As a result, we have
\begin{align*}
\widehat{\BS}_i^{\eig}= \argmin_{\BS_i^{\top}\BS_i = \BPi_i}-\lag \BS_i, \widehat{\BS}_i^{\eig}\rag, \qquad~~ \widehat{\BS}_i^{\eig} = \argmin_{\BS_i^{\top}\BS_i = \BPi_i} \lag \BA_i, \BS_i\rag
\end{align*}
Again, by using $\BPhi  = \BS_i\BPi_i^{-1/2}$, then
\[
\min_{\BS_i^{\top}\BS_i = \BPi_i} \lag \BA_i, \BS_i\rag = \min_{\BPhi \in\Od(d)} \lag \BA_i\BPi_i^{1/2}, \BPhi \rag
\Longrightarrow
\argmin_{\BPhi \in\Od(d)} \lag \BA_i\BPi_i^{1/2}, \BPhi \rag = \widehat{\BS}_i^{\eig}\BPi_i^{-1/2} = \PP(\widehat{\BS}_i^{\eig}).
\]
This finishes the proof.
\end{proof}

The derivation for the distribution of spectral estimator is quite similar to its MLE counterpart.
Define
$f(\BPhi) : = \lag \BA, \BPhi\rag$ such that $\BPhi\in\Od(d)$ and
\begin{equation}\label{def:A_spectral}
\BA :  = \BA_i \BPi_i^{1/2}= -\sum_{j=1}^n \BC_{ij}\widehat{\BS}_j^{\eig} \BPi_i^{1/2}= - ( \BZ_i\BZ^{\top}\BSeig + \sigma\BW_i^{\top}\BSeig) \BPi_i^{1/2}.
\end{equation}
Lemma~\ref{eq:spr_local} implies that the global minimizer to $f(\BPhi)$ is $\widehat{\BZ}_i^{\eig}$. 
We will show that the one-step Newton iteration applied to $f(\cdot)$ at $\BZ_i$ is asymptotically Gaussian and approximates $\widehat{\BZ}^{\eig}_{i}$ up to an error of order $O \lp \SNR^{-2}(\sqrt{d} + \sqrt{\log n})^2 \rp$. 
Based on Proposition~\ref{prop:2nderror}, the one-step Newton iteration at $\BZ_{i}$ is given by $\widetilde{\BZ}^{\eig}_{i}:= R_{\BZ_i}(\widetilde{\BV}^{\rm eig}_i) \in \Od(d)$ where
\begin{equation}\label{def:ntsp}
\begin{aligned}
& \widetilde{\BV}^{\rm eig}_{i} := -\BH_i^{{\rm eig}-1}(\BA_{i}\BPi^{1/2}_{i} + \BH^{\rm eig}_i \BZ_i), \\ 
& \BH_i^{{\rm eig}} := - \frac{1}{2} (\BA\BZ_i^{\top} + \BZ_i\BA^{\top})= -\frac{1}{2} (\BA_{i}\BPi^{1/2}_{i}\BZ_i^{\top} + \BZ_i\BPi_i^{1/2}\BA^\top_{i}),
\end{aligned}
\end{equation}
and $\widetilde{\BV}_i^{\eig}$ is on the tangent space of $\Od(d)$ at $\BZ_i.$

Without loss of generality, we choose $\BSeig$ such that $\PP(\BZ^{\top}\BSeig) = \I_d$ since for any $\BQ\in\Od(d)$, it holds $\widehat{\BS}^{\eig}\BQ (\widehat{\BS}^{\eig}\BQ)^{\top} = \widehat{\BS}^{\eig}\widehat{\BS}^{\eig\top}$. As a result, we have  $d_F(\widehat{\BS}^{\eig}, \BZ) = \|\widehat{\BS}^{\eig} - \BZ\|_F$ where $\widehat{\BS}^{\eig}$ is the global minimizer to~\eqref{eq:spectral}.

\subsection{Supporting lemmas}
\begin{lemma}{\bf~\cite[Theorem 5.1]{L22b}}\label{lem:pi}
Under the conditions of Theorem~\ref{thm:spectral_main}, it holds with high probability that
\[
|\sigma_{j}(\BSeig_i) - 1| \lesssim \sigma n^{-1/2}(\sqrt{d} + \sqrt{\log n}) \lesssim \SNR^{-1}d^{-1/2} (\sqrt{d} + \sqrt{\log n}),~~\forall 1\leq j\leq d,
\]
where $\sigma_j(\cdot)$ denotes the $j$-th largest singular value.
Thus  $\BPi_i = \widehat{\BS}_i^{\eig\top}\widehat{\BS}_i^{\eig}$ satisfies
\[
\norm{\BPi_{i} - \BI_{d}} \lesssim \SNR^{-1} d^{-1/2} (\sqrt{d} + \sqrt{\log n}),~~~\norm{\BPi_{i}^{1/2} - \BI_{d}} \lesssim \SNR^{-1} d^{-1/2} (\sqrt{d} + \sqrt{\log n}).
\] 
\end{lemma}
Lemma~\ref{lem:pi} implies an estimation on the smallest and largest eigenvalues $\pi_{\min}$ $(\pi_{\max})$ of $\BPi_i.$
The condition number  of $\BPi_i$ satisfies
\begin{equation}\label{def:kappa}
\kappa : =\frac{\pi_{\max}}{\pi_{\min}}  \leq \frac{1 + c \SNR^{-1}d^{-1/2} (\sqrt{d} + \sqrt{\log n})}{1 - c \SNR^{-1}d^{-1/2} (\sqrt{d} + \sqrt{\log n})}
\end{equation}
for some $c>0$ where $\pi_{\max} = \sigma_1^2(\widehat{\BS}_i^{\eig})$ and $\pi_{\min} = \sigma_d^2(\widehat{\BS}_i^{\eig}).$

\begin{lemma}\label{lem:assumption2}
Under the conditions of Theorem~\ref{thm:spectral_main},   the spectral estimator $\widehat{\BS}^{\eig}$ and its rounding $\widehat{\BZ}^{\eig}$ satisfy
\begin{equation}\label{eq:assumption2}
\begin{aligned}
& \min_{\BQ\in\Od(d)}\max_{1\leq i\leq n}\| \widehat{\BZ}_i^{\eig} - \BZ_i \BQ \|_F\lesssim \sigma n^{-1/2}d = \SNR^{-1}\sqrt{d}, \\
& \min_{\BQ \in \Od(d)} \|\widehat{\BS}^{\eig}- \BZ\BQ\|_F\leq C\cdot \SNR^{-1}\sqrt{nd}, \\
& \|(\BZ^{\top}\BSeig - n\BQeig) \|_F \leq \frac{1}{2} \|\BSeig - \BZ\BQeig\|_F^2 \lesssim   n\SNR^{-2}(\sqrt{d}+\sqrt{\log n})^2, 
\end{aligned}
\end{equation}
with high probability where $\BQ^{\eig} = \PP(\BZ^{\top}\widehat{\BS}^{\eig})$.
\end{lemma}
The first one follows from~\cite[Theorem 3.1]{L22b}: 
\begin{equation}
\min_{\BQ\in\Od(d)}\max_{1\leq i\leq n}\| \widehat{\BZ}_i^{\eig} - \BZ_i \BQ \|_F\lesssim \sigma \sqrt{\frac{d}{n}}(\sqrt{d} +\sqrt{\log n}) = \SNR^{-1}(\sqrt{d} +\sqrt{\log n}).
\end{equation}
The second equation in~\eqref{eq:assumption2} follows from Davis-Kahan theorem.
By applying Lemma~\ref{lem:zs} and~\ref{lem:pi}, we get the estimation of $\| \BZ^{\top}\BSeig - n\BQeig \|_F$ in~\eqref{eq:assumption2}.

\begin{lemma}\label{lem:ws_spectral}
Under the conditions of Theorem~\ref{thm:main}, $\BW_i^{\top}\BSeig$ satisfies 
\begin{align*}
\BW_i^{\top}\BSeig \overset{d}{=} \sqrt{n-1} \BG_i + O\lp \sqrt{nd}  \SNR^{-1}(\sqrt{d}+\sqrt{\log n}) \rp,~~1\leq i\leq n, \\
\BW_i^{\top}\BSeig\BPi_{i}^{1/2} \overset{d}{=} \sqrt{n-1} \BG_i + O\lp \sqrt{nd}  \SNR^{-1}(\sqrt{d}+\sqrt{\log n})^2 \rp,~~1\leq i\leq n,
\end{align*}
where $\BG_i$ is a $d\times d$ random matrix with i.i.d. standard Gaussian entries. As a result, it holds $\max_{1\leq i\leq n}\| \BW_i^{\top}\widehat{\BS}^{\eig} \|_F \leq 3 \sqrt{n}(\sqrt{d} + \sqrt{\log n})$
with high probability.
\end{lemma}

\begin{proof}[\bf Proof of Lemma~\ref{lem:ws_spectral}]
Due to the statistical dependency between $\BW_i$ and $\BSeig$, we let $\BS^{\eig(i)}$ be the spectral estimator associated with $\BC^{(i)}$ where $\BC^{(i)}$ is equal to $\BC$ except with the noise on the $i$-th block row and column of $\BC$ removed, as shown in~\eqref{def:Ci}.
From~\cite[Equation (5.8)]{L22b}, we know that the distance between $\BS^{\eig}$ and $\BS^{\eig(i)}$ is very small:
\[
\norm{\BSeig - \BSeigloo\BQ^{i}}_F \lesssim \sigma n^{-1/2}(\sqrt{d} + \sqrt{\log n}) \sqrt{d}\max_{1\leq i\leq n} \|\widehat{\BS}_i^{\eig}\| 
\lesssim \SNR^{-1}(\sqrt{d} + \sqrt{\log n}) 
\]
where $\BQ^i = \PP(\widehat{\BS}^{\eig(i)\top} \widehat{\BS}^{\eig})$ and $\|\widehat{\BS}_i^{\eig}\| \leq 2$ follows from Lemma~\ref{lem:pi}.
We proceed with the estimation of $\BW_i^{\top}\widehat{\BS}^{\eig}$ which is decomposed into
\begin{align*}
\BW_i^{\top}\BSeig 
& = \BW_i^{\top}(\BSeig - \BSeigloo\BQ^i)  + \BW_i^{\top} \BSeigloo(\BQ^i - \I_d) + \BW_i^{\top} \BSeigloo  = : T_1' + T_2' + T_3'
\end{align*}
where $\BQ^i = \PP(\BSeigloot\BSeig )$.
For the first term, we note that
\[
\|T_1'\|_F \le \norm{\BW_i} \cdot \norm{\BSeig - \BSeigloo\BQ^i}_F \lesssim \sqrt{nd} \cdot\SNR^{-1}(\sqrt{d}+\sqrt{\log n})
\]
where $\|\BW_i\|\leq 3\sqrt{nd}$ with high probability.
For the second term, it holds that
\begin{align*}
\|T_2'\|_F & \leq \|\BW_i\|\|\BSeig\|\|\BQ^i-\I_d\|_F \lesssim \sqrt{nd} \cdot \sqrt{n}\cdot n^{-1/2} \SNR^{-1}(\sqrt{d}+\sqrt{\log n}) \\
& \lesssim \sqrt{nd} \SNR^{-1}(\sqrt{d}+\sqrt{\log n}) 
\end{align*}
where 
$\|\BQ^i-\I_d\|_F \lesssim n^{-1/2} \SNR^{-1}(\sqrt{d}+\sqrt{\log n})$
follows the exact same argument as~\eqref{eq:QminusId}.
For the third term, we have
\[
T_3' = \BW_i^{\top} \BSeigloo = \BW_{i}^{\top}\BZ\BQ^{(i)} + \BW_{i}^{\top}(\BSeigloo - \BZ\BQ^{(i)})
\]
where $\BQ^{(i)} := \PP(\BZ^{\top}\widehat{\BS}^{\eig (i)})$.
Unlike $\widehat{\BS}$, each block of $\widehat{\BS}^{\eig(i)}$ is non-orthogonal and $\BW_{ii} = 0.$
 Note that $\BW_{i}^{\top}\BZ\BQ^{i}$ is the sum of $n-1$ independent standard Gaussian matrices and then we have
$\BW_i^{\top}\BZ\BQ^{(i)} \overset{d}{= }\sqrt{n-1}\BG_i$
where $\BQ^{(i)}$ is independent of $\BW_i$ and $\BG_i$ is a $d\times d$ Gaussian random matrix with i.i.d. entries.
Lemma 5.8 in~\cite{L22b} gives
\[
\norm{\BW_{i}^{\top}(\BSeigloo - \BZ\BQ^{(i)})}_{F} \lesssim  \norm{\BSeigloo - \BZ\BQ^{(i)}}_F (\sqrt{d}+\sqrt{\log n}) \lesssim \sqrt{nd}\SNR^{-1}(\sqrt{d}+\sqrt{\log n}).
\]
This is because $\BW_i$ and $\BSeigloo - \BZ\BQ^{(i)}$ are independent, and 
the squared Frobenius norm of $\BW_{i}^{\top}(\BSeigloo - \BZ\BQ^{(i)})$ is controlled by a $\chi_d^2$-random variable times $\norm{\BSeigloo - \BZ\BQ^{(i)}}_F^2$, and  $\norm{\BSeigloo - \BZ\BQ^{(i)}}_F \lesssim \sqrt{nd}\SNR^{-1}.$
Therefore, we have
\[
T_3'  = \sqrt{n-1}\BG_i + O(\sqrt{nd}\SNR^{-1}(\sqrt{d} + \sqrt{\log n}))
\]
Combining $T_1'$, $T_2'$ with $T_3'$ immediately implies
\begin{equation} \label{eq:ws_spectral}
\begin{aligned}
\BW_i^{\top}\widehat{\BS}^{\eig} & \overset{d}{=} \sqrt{n-1}\BG_i + O(\sqrt{nd}\SNR^{-1}(\sqrt{d} + \sqrt{\log n})), \\
    \norm{\BW_i^{\top}\BSeig}_{F} & \le \norm{T'_1}_F + \norm{T'_2}_F + \norm{T'_3}_F\lesssim \sqrt{nd}(\sqrt{d} +\sqrt{\log n})
\end{aligned}
\end{equation}
where $\SNR$ is greater than $\sqrt{d} + \sqrt{\log n}$ and $\|\BG_i\|_F\lesssim \sqrt{d} + \sqrt{\log n}$ holds with high probability.
For $\BW_i^{\top}\widehat{\BS}^{\eig}\BPi_i$, it take only one more step:
\[
\BW_i^{\top}\widehat{\BS}^{\eig}\BPi_i = \BW_i^{\top}\widehat{\BS}^{\eig}+ \underbrace{\BW_i^{\top}\widehat{\BS}^{\eig}(\BPi_i-\I_d)}_{T_4'}.
\]
The fourth term is upper bounded via 
\[
\begin{aligned}
	\norm{T'_4}_{F} &= \norm{\BW_i^{\top}\BSeig(\BPi_{i}^{1/2}-\BI_{d})}_{F} \le \norm{\BW_i^{\top}\BSeig}_{F}\norm{\BPi_{i}^{1/2}-\BI_{d}}\\
&\lesssim \sqrt{n}d \cdot  \SNR^{-1} d^{-1/2} (\sqrt{d} + \sqrt{\log n}) = \sqrt{n}\SNR^{-1}(\sqrt{d} + \sqrt{\log n})^2.
\end{aligned}
\]
This gives the estimation of $\BW_i^{\top}\widehat{\BS}^{\eig}\BPi_i.$
\end{proof}

\begin{lemma}\label{lem:Hesseig}
Under the conditions of Theorem~\ref{thm:spectral_main}, it holds with high probability that
\[
\begin{aligned}
	\|\BH_i^{\eig} - n\I_d\| & \lesssim  \sqrt{\pi_{\max}}n\SNR^{-1}(\sqrt{d} + \sqrt{\log n})
\end{aligned}
\] 
where $\BH_i^{\eig}$ corresponds to the Hessian of $f(\BPhi) = \lag \BA,\BPhi\rag$ in with $\BA$ in~\eqref{def:A_spectral}.
\end{lemma}
\begin{proof}[\bf Proof of Lemma~\ref{lem:Hesseig}]
Without loss of generality, we choose $\BSeig$ such that it satisfies $\BQeig = \I_d$, and then $\| \BZ^{\top} \BSeig - n\I_d\| \lesssim n \SNR^{-2}(\sqrt{d} + \sqrt{\log n})^2$
which follows from Lemma~\ref{lem:assumption2}.
We have
\[
\BH^{\rm eig}_i = \frac{1}{2}  \BZ_i (\BZ^{\top}\BSeig\BPi^{1/2}_{i} + \BPi^{1/2}_{i}\BSeigt\BZ)\BZ_i^{\top} + \frac{\sigma}{2} \lp \BW_i^{\top}\BSeig\BPi^{1/2}_{i} \BZ_i^{\top} + \BZ_i \BPi^{1/2}_{i}\BSeigt\BW_i \rp
\]
which follows from~\eqref{def:Aisp} and~\eqref{def:ntsp}.
Then  for $\BH_i^{\eig}$, it holds that
\begin{align*}
& \| \BH_i^{\eig} - n\BPi_i^{1/2} \| \\
& ~~\leq \frac{1}{2} \left\| \BZ_i ((\BZ^{\top}\BSeig - n\I_d)\BPi^{1/2}_{i} + \BPi^{1/2}_{i}(\BSeigt\BZ  - n\I_d))\BZ_i^{\top} \right\| \\
&~~\qquad +\frac{\sigma}{2} \left\| \BW_i^{\top}\BSeig\BPi^{1/2}_{i} \BZ_i^{\top} + \BZ_i \BPi^{1/2}_{i}\BSeigt\BW_i\right\| \\
&~~ \leq \left(\| \BZ^{\top}\BSeig - n\I_d  \|_F + \sigma \| \BW_i^{\top}\widehat{\BS}^{\eig} \|_F\right) \|\BPi_i^{1/2}\|  \lesssim \sqrt{\pi_{\max}}
n\SNR^{-1} (\sqrt{d} +\sqrt{\log n})
\end{align*}
where Lemma~\ref{lem:ws_spectral} gives $\norm{\BW_i^{\top}\BSeig}_{F} \lesssim \sqrt{nd}(\sqrt{d}+\sqrt{\log n})$.
\end{proof}

\subsection{Proof of Proposition~\ref{prop:tildeV_spectral}}
\begin{proof}[\bf Proof of Proposition~\ref{prop:tildeV_spectral}]
Note that the Newton step is given by
\begin{align*}
-\widetilde{\BV}_i^{\eig}  & =  \BH_i^{{\rm eig}-1}(\BA_{i}\BPi^{1/2}_{i} + \BH^{\rm eig}_i \BZ_i)  = \frac{1}{2}\BH^{{\rm eig} -1}_i ( \BA_{i}\BPi^{1/2}_{i}\BZ_i^{\top} - \BZ_i\BPi_i^{1/2}\BA^\top_{i} )\BZ_i \\
& = \frac{1}{2}\left(\BH^{{\rm eig} -1}_i  - n^{-1}\I_d\right) ( \BA_{i}\BPi^{1/2}_{i}\BZ_i^{\top} - \BZ_i\BPi_i^{1/2}\BA^\top_{i} )\BZ_i  + \frac{1}{2n} ( \BA_{i}\BPi^{1/2}_{i}\BZ_i^{\top} - \BZ_i\BPi_i^{1/2}\BA^\top_{i})\BZ_i 
\end{align*}
where $\BH_i^{{\rm eig}} = -\frac{1}{2} (\BA_{i}\BPi^{1/2}_{i}\BZ_i^{\top} + \BZ_i\BPi_i^{1/2}\BA^\top_{i})$.

We will estimate each term above. For the first term in $\widetilde{\BV}_i^{\eig}$, note that
\begin{equation} \label{eq: T_1_1_spectral}
\begin{aligned}
 &\| \BH^{{\rm eig} -1}_i  - n^{-1}\I_d\| = n^{-1} \|\BH^{{\rm eig} -1}_i\| \|\BH^{{\rm eig} }_i - n\I_d\|\\
&~~~~~~ \lesssim \frac{1}{n^2\sqrt{\pi_{\min}}}\cdot n\sqrt{\pi_{\max}} \SNR^{-1}(\sqrt{d} + \sqrt{\log n}) \lesssim \sqrt{\kappa}n^{-1} \SNR^{-1}(\sqrt{d} +\sqrt{\log n})
\end{aligned}
\end{equation}
where $\lambda_{\min}(\BH_i^{\eig}) \geq \sqrt{\pi_{\min}}n$ and 
the estimation of $\|\BH_i^{\eig} - n\I_d\|$ are given by Lemma~\ref{lem:Hesseig} under $\SNR\gtrsim \sqrt{\kappa} (\sqrt{d} + \sqrt{\log n})$.

Let
\begin{align*}
 &\BA_{i}\BPi^{1/2}_{i}\BZ_i^{\top} - \BZ_i\BPi_i^{1/2}\BA^\top_{i}  := T_1 + T_2 \\
&\qquad = \underbrace{ \BZ_i(\BZ^{\top}\BSeig\BPi^{1/2}_{i} - \BPi^{1/2}_{i}\BSeigt\BZ)\BZ_i^{\top} }_{T_1} + \underbrace{\sigma (\BW_i^{\top}\BSeig\BPi^{1/2}_{i} \BZ_i^{\top} - \BZ_i \BPi^{1/2}_{i}\BSeigt\BW_i) }_{T_2}
\end{align*}
and then
\[
-\widetilde{\BV}_i^{\eig} = \frac{1}{2}(\BH_i^{\eig-1} - n^{-1}\I_d)(T_1 + T_2)\BZ_i + \frac{1}{2n}(T_1+T_2)\BZ_i.
\]
It suffices to estimate $T_1$ and $T_2$ respectively.
For $T_1$,  Lemma~\ref{lem:assumption2} gives 
\begin{align*}
\|T_1\|_F & \leq \left\| \BZ^{\top}\BSeig\BPi^{1/2}_{i} - \BPi^{1/2}_{i} \BSeigt \BZ  \right\|_F  = \left\| (\BZ^{\top}\BSeig - n\I_d)\BPi^{1/2}_{i} - \BPi^{1/2}_{i} (\BSeigt \BZ -n\I_d) \right\|_F  \\
& \leq 2 \sqrt{\pi_{\max}} \| \BZ^{\top}\widehat{\BS}^{\eig} -n\I_d\|_F \leq  2n \sqrt{\pi_{\max}} (\SNR^{-1}(\sqrt{d}+\sqrt{\log n}))^2.
\end{align*}
For $T_2$, Lemma~\ref{lem:ws_spectral} implies
\[
(\BW_i^{\top}\BSeig\BPi^{1/2}_{i} \BZ_i^{\top} - \BZ_i \BPi^{1/2}_{i}\BSeigt\BW_i) \BZ_i = \sqrt{n-1}(\BG_i - \BG_i^{\top}) \BZ_i + O\lp \sqrt{nd} \SNR^{-1}(\sqrt{d} + \sqrt{\log n})^2\rp
\]
where $\BG_i$ is a $d\times d$ symmetric Gaussian random matrix. Then 
\begin{align*}
 T_2 & \overset{d}{=} \sqrt{n-1}\sigma (\BG_i-\BG_i^{\top})\BZ_i + O(\sigma\sqrt{nd}\SNR^{-1}(\sqrt{d} +\sqrt{\log n})^2),\\
\|T_2\|_F & \lesssim \sigma \sqrt{nd} (\sqrt{d} +\sqrt{\log n}) \lesssim n \SNR^{-1}(\sqrt{d} +\sqrt{\log n}).
\end{align*}
As a result, we have
\begin{align*}
T_1 + T_2 &  \overset{d}{=}  \sqrt{n-1}\sigma (\BG_i-\BG_i^{\top})\BZ_i + O(\sigma\sqrt{nd}\SNR^{-1}(\sqrt{d} +\sqrt{\log n})^2), \\
\|T_1 + T_2\|_F &  = \| \BA_{i}\BPi^{1/2}_{i}\BZ_i^{\top} - \BZ_i\BPi_i^{1/2}\BA^\top_{i}\|_F \lesssim n\SNR^{-1}(\sqrt{d} +\sqrt{\log n})
\end{align*}
provided that $\SNR\gtrsim \sqrt{\kappa}(\sqrt{d} +\sqrt{\log n}).$

Finally, the first term in $\widetilde{\BV}_i^{\eig}$ is bounded by
\begin{align*}
& \left\| \left(\BH^{{\rm eig} -1}_i  - n^{-1}\I_d\right) ( T_1 + T_2)\BZ_i\right\|_F \leq \left\| \BH^{{\rm eig} -1}_i  - n^{-1}\I_d\right\| \cdot \| T_1 + T_2\|_F \\
& ~~~ \lesssim \sqrt{\kappa}n^{-1} \SNR^{-1}(\sqrt{d} +\sqrt{\log n}) \cdot n\SNR^{-1}(\sqrt{d} +\sqrt{\log n}) \lesssim \sqrt{\kappa} (\SNR^{-1}(\sqrt{d} +\sqrt{\log n}))^2
\end{align*}
and 
\begin{align*}
\widetilde{\BV}_i^{\eig} & = \frac{1}{2} (\BH_i^{\eig -1} - n^{-1}\I_d) (T_1 + T_2)\BZ_i + \frac{1}{2n} (T_1 + T_2)\BZ_i \\
& \overset{d}{=} \frac{\sqrt{n-1}\sigma}{2n}(\BG_i - \BG_i^{\top}) \BZ_i + O(\sqrt{\kappa}  (\SNR^{-1}(\sqrt{d} +\sqrt{\log n}))^2)
\end{align*}
where the leading term is given by $T_2/n.$
As a result, the norm of $\widetilde{\BV}_i^{\eig}$ is bounded by
\[
\|\widetilde{\BV}_i^{\eig}\|_F \lesssim \sigma\sqrt{\frac{d}{n}}(\sqrt{d} +\sqrt{\log n}) \lesssim \SNR^{-1}(\sqrt{d} +\sqrt{\log n})
\]
under $\SNR\gtrsim\sqrt{\kappa}(\sqrt{d} +\sqrt{\log n}).$
\end{proof}

\subsection{Proof of Proposition~\ref{prop:newtonstep_spectral}}
\begin{proof}[\bf Proof of Proposition~\ref{prop:newtonstep_spectral}]
Consider $f(\BPhi) = \lag \BA, \BPhi\rag$ where $\BPhi\in\Od(d)$ and $\BA = \BA_i\BPi_i^{1/2}$ is in~\eqref{def:A_spectral}.
Let $\widehat{\BV}^{\rm eig}_i$ and $\widetilde{\BV}^{\rm eig}_{i} $ be the two tangent vectors such that 
\[
R_{\BZ_i}(\widehat{\BV}^{\rm eig}_i) = \widehat{\BZ}^{\eig}_i,~~~~ R_{\BZ_i}(\widetilde{\BV}^{\rm eig}_i) = \widetilde{\BZ}^{\eig}_i
\] 
where $\widetilde{\BV}_i^{\eig}$ is given in~\eqref{def:ntsp}.
Lemma~\ref{lemma: Lipschitz} guarantees the existence  of $\widehat{\BV}_i^{\eig}$ if $\|\widehat{\BZ}^{\eig}_i - \BZ_i\| < 1$ which follows from
\begin{equation}\label{eq:Zeig_diff}
\|\widehat{\BZ}^{\eig}_i - \BZ_i\| \leq \|\widehat{\BZ}^{\eig}_i - \BZ_i\|_F \lesssim \SNR^{-1} (\sqrt{d} +\sqrt{\log n}) < \frac{1}{8}
\end{equation}
if $\SNR\gtrsim \sqrt{d} +\sqrt{\log n}.$

Similar to the proof in Section~\ref{ss:propkey}, the idea of controlling $\|\widehat{\BZ}_i^{\eig} - \widetilde{\BZ}_i^{\eig}\|_F$ is to apply Lemma~\ref{prop:2nderror}, i.e.,
\begin{equation}
    \begin{aligned} \label{eqn: prop_spectral_main}
	\|\widehat{\BZ}^{\eig}_{i} - \widetilde{\BZ}^{\eig}_{i}\|_F 
	& \leq \frac{3}{2} \left( \frac{\|\BA_{i}\BPi^{1/2}_{i}\|}{\lambda_{\min}(\BH^{\rm eig}_i)} +  \sqrt{  \frac{\|\BA_{i}\BPi^{1/2}_{i}\|}{\lambda_{\min}(\BH^{\rm eig}_i)}} \right) (\|\widehat{\BV}^{\rm eig}_i\|_F^2  + \|\widetilde{\BV}^{\rm eig}_i\|_F^2 )
\end{aligned}
\end{equation}
where $\BA = \BA_i\BPi_i^{1/2}$.
Therefore, it suffices to bound $\norm{\BA_{i}\BPi^{1/2}_{i}}$, $\|\widehat{\BV}^{\rm eig}_{i}\|_F$, and $\|\widetilde{\BV}^{\rm eig}_{i}\|_F.$ Using Lemma~\ref{lem:pi}, we have that
\begin{equation}\label{eqn: spectral norm of A spectral}
\begin{aligned}
	\norm{\BA_{i}\BPi^{1/2}_{i}} & \le \norm{\BPi_{i}}^{1/2}\norm{\sum_{j=1}^n\BC_{ij}\BSeig_{j}} \leq \sqrt{\pi_{\max}} \left\|  \BZ_i\BZ^{\top}\BSeig + \sigma\BW_i^{\top}\BSeig \right\| \\
		&\lesssim \sqrt{\pi_{\max}} \left( n +\sigma \|\BW_i^{\top}\BSeig \|_F \right) \lesssim \sqrt{\pi_{\max}} (n + 3\sigma\sqrt{nd}(\sqrt{d} +\sqrt{\log n})) \lesssim \sqrt{\pi_{\max} }n
\end{aligned}	
\end{equation}
where the estimation of $ \|\BW_i^{\top}\BSeig \|_F \lesssim \sqrt{nd}(\sqrt{d} +\sqrt{\log n})$ is given in~\eqref{lem:ws_spectral} and $\SNR \gtrsim \sqrt{d} +\sqrt{\log n}.$

Using the Lipschitz continuity of $R_{\BZ_{i}}^{-1}$ (Lemma~\ref{lemma: Lipschitz}) with~\eqref{eq:Zeig_diff}, $\|\widehat{\BV}^{\rm eig}_{i}\|_F$ is upper bounded via 
\begin{align} \label{eqn: F_norm of V_hat spectral}
	\|\widehat{\BV}^{\rm eig}_{i}\|_F \le 2 \norm{\widehat{\BZ}_{i}^{\eig} - \BZ_{i}}_{F}  \lesssim  \SNR^{-1}(\sqrt{d}+\sqrt{\log n}).
\end{align}
Proposition~\ref{prop:tildeV_spectral} gives 
$\|\widetilde{\BV}_i\|_F \lesssim \SNR^{-1}(\sqrt{d} +\sqrt{\log n}).$
Finally, we have
\[
\max\{\|\widehat{\BV}_i\|_F, \|\widetilde{\BV}_i\|_F\} \lesssim \SNR^{-1}(\sqrt{d} +\sqrt{\log n})
\]
and combining~\eqref{eqn: spectral norm of A spectral} with $\lambda_{\min}(\BH_i^{\eig})\geq  \sqrt{\pi_{\min}}n/2$ gives
$\|\BA_i\BPi_i^{1/2}\|/\lambda_{\min}(\BH_i^{\eig})  \lesssim \sqrt{\kappa}.$
Applying these estimations~\eqref{eqn: prop_spectral_main} leads to
$\norm{\widehat{\BZ}_{i} - \widetilde{\BZ}_{i}}_{F} \leq \|\widehat{\BV}_{i} - \widetilde{\BV}_{i} \|\lesssim \kappa^{1/2}  \SNR^{-2}(\sqrt{d}+\sqrt{\log n})^2.$
\end{proof}

{
\section{Proof of Corollary~\ref{cor: CR mle}}


\begin{proposition} \label{prop: estimation error of sigma2}
	For $\widehat{\BS} \in \{\widehat{\BZ},\BZeig\}$, the variance estimator $\widehat{\sigma}^2$ obtained from~\eqref{eq: variance estimation} satisfies
	\[
	\begin{aligned}
		\abs{\frac{\widehat{\sigma}^2}{\sigma^2} - 1} \lesssim \frac{1}{n}\left( 1 + \frac{\sqrt{\log n}}{d}\right).
	\end{aligned}
	\] 
\end{proposition}

\begin{proof}[\bf Proof of Proposition~\ref{prop: estimation error of sigma2}]
The estimation error of $\widehat{\sigma}^2$ is given by
\begin{equation} \label{eq: estimation error of sigma2}
	\begin{aligned}
	& \abs{ \widehat{\sigma}^2 - \sigma^2}  = \abs{ \frac{1}{n(n-1)d^2} \|\BC - \widehat{\BS}\widehat{\BS}^{\top}\|_F^2 - \sigma^2 } \\
	&\le \sigma^2\abs{\frac{\|\BW\|_F^2}{n(n-1)d^2} - 1}  
	+ \frac{1}{n(n-1)d^2} \left(  2\sigma |\lag \BZ\BZ^{\top} - \widehat{\BS}\widehat{\BS}^{\top}, \BW\rag| + \| \BZ\BZ^{\top} - \widehat{\BS}\widehat{\BS}^{\top}\|_F^2 \right) \\
	\end{aligned}
\end{equation}
where
\[
\|\BC - \widehat{\BS}\widehat{\BS}^{\top}\|_F^2 = \sigma^2 \|\BW\|_F^2 + 2\sigma \lag \BZ\BZ^{\top} - \widehat{\BS}\widehat{\BS}^{\top}, \BW\rag + \|\BZ\BZ^{\top} - \widehat{\BS}\widehat{\BS}^{\top}\|_F^2
\]

Note that $2^{-1}\|\BW\|_F^2$ is $\chi^2$-distribution of degree $n(n-1)d^2/2$ since the diagonal blocks are zero. Then
\[
\Pr\left(  \frac{1}{2}\left|  \|\BW\|_F^2 - n(n-1)d^2  \right| \geq t  \right) \leq 2\exp\left( -\frac{t^2}{4n(n-1)d^2}\right) \vee \exp\left( -\frac{t}{8}\right).
\]
Then by picking $t = 2\sqrt{n(n-1)d^2 \log n} $, then 
\begin{align*}
& \left| \|\BW\|_F^2 - n(n-1)d^2 \right| \leq 4\sqrt{n(n-1)d^2 \log n} \\
&~~~ \Longleftrightarrow \sigma^2 \left| \frac{\|\BW\|_F^2}{n(n-1)d^2} - 1 \right| \leq 4\sigma^2 \sqrt{\frac{\log n}{n(n-1)d^2} } = 4\SNR^{-2}d^{-2}\sqrt{\frac{n\log n}{n-1}}
\end{align*}
with high probability.


Note that for both MLE (Proposition 5.6 in \cite{L23b}) and spectral estimator (Lemma~\ref{lem:assumption2}) $\widehat{\BS} \in \{ \widehat{\BZ}, \widehat{\BZ}^{\eig}\}$, it holds 
\[
\norm{\widehat{\BS} - \BZ\BQ}_{F} \lesssim \SNR^{-1}\sqrt{nd}.
\]
Then the second and third term are upper bounded via 
\[
\begin{aligned}
	\norm{\BZ\BZ^{\top} - \widehat{\BS}\widehat{\BS}^{\top}}_{F} & = \norm{\BZ\BQ(\BZ\BQ)^{\top} - \widehat{\BS}(\BZ\BQ)^{\top} + \widehat{\BS}(\BZ\BQ)^{\top} - \widehat{\BS}\widehat{\BS}^{\top}}_{F}\\
	&\le \norm{\BZ}\norm{\widehat{\BS}-\BZ\BQ}_F + \norm{\widehat{\BS}}\norm{\widehat{\BS}-\BZ\BQ}_F\\
	&\lesssim 2\sqrt{n} \cdot \SNR^{-1}\sqrt{nd}\\
	&\lesssim \SNR^{-1} n\sqrt{d}
\end{aligned}
\]
and
\[
\begin{aligned}
\sigma\left \langle\BW, \BZ\BZ^{\top} - \widehat{\BS}\widehat{\BS}^{\top} \right \rangle  &= \sigma\left \langle\BW,  \lp \BZ\BQ - \widehat{\BS} \rp\lp \BZ\BQ + \widehat{\BS} \rp^{\top} \right \rangle\\ 
&= \sigma\left \langle\BW \lp \BZ\BQ + \widehat{\BS} \rp,   \BZ\BQ - \widehat{\BS}  \right \rangle\\
	&\le \sigma \norm{\BW \lp  \widehat{\BS} + \BZ\BQ \rp}_{F}\cdot \norm{\widehat{\BS} - \BZ\BQ}_{F}\\
	&\lesssim \SNR^{-1}\sqrt{\frac{n}{d}} \cdot 3\sqrt{nd} \cdot 2\sqrt{nd} \cdot \SNR^{-1}\sqrt{nd}\\
	&\lesssim \SNR^{-2} n^2d
\end{aligned}
\]
Thus, \eqref{eq: estimation error of sigma2} implies that the relative error between $\widehat{\sigma}^2$ and $\sigma^2$ is given by
\begin{align}\label{eq: relative error of sigma2}
	\abs{\frac{\widehat{\sigma}^2 - \sigma^2}{\sigma^2}} \lesssim  \frac{\SNR^{-2}d^{-1} + \SNR^{-2}d^{-2}\sqrt{\log n} }{\SNR^{-2} d^{-1}n} = \frac{1}{n}\left( 1 + \frac{\sqrt{\log n}}{d}\right).
\end{align}
\end{proof}
}

\bibliographystyle{abbrv}

\end{document}